\documentclass[a4paper,12pt]{amsart}
\usepackage[letterpaper, margin=1.2in]{geometry}
\usepackage{latexsym}
\usepackage{amssymb}
\usepackage{amsmath}
\usepackage{amsfonts}
\usepackage[table]{xcolor}
\usepackage{pgfplots}
\usepackage{color}
\usepackage{multirow}
\usepackage{graphicx}
\usepackage{txfonts,pxfonts} 
\usepackage{tikz}
\usetikzlibrary{calc,arrows}
\usetikzlibrary{calc}
\usepackage{verbatim}
\everymath{\displaystyle}
\newtheorem{thm}{Theorem}[section]

\newtheorem{lem}[thm]{Lemma}

\newtheorem{cor}[thm]{Corollary}

\theoremstyle{definition}

\newcommand{\Z}{\mathbb Z}

\newcommand{\Q}{\mathbb Q}

\newcommand{\F}{\mathbb F}

\def\F{\mathbb{F}}

\def\ol{\overline}

\def\la{\lambda}

\def\om{\omega}

\def\md#1{\ \mbox{\rm(mod }{#1})}

\def\ph{\phi}

\newcounter{cs}
\stepcounter{cs}
\newcommand{\casos}{\begin{itemize}}
	\newcommand{\fcasos}{\end{itemize}\setcounter{cs}{1}}

\newfont{\tit}{cmr12 scaled \magstep3}
\title[On index divisors and monogenity of cetain number fields defined by $x^{12}+ax^m+b$]{ On index divisors and monogenity of certain number fields defined by $x^{12}+ax^m+b$}
\author[L. El Fadil]{Lhoussain El Fadil}
\author[O. Kchit]{Omar Kchit}
\address{Faculty of Sciences Dhar El Mahraz, P.O. Box  1874 Atlas-Fes, Sidi Mohamed ben Abdellah University,  Morocco}
\email{lhoussain.elfadil@usmba.ac.ma, \,\, orcid: 0000.0003.4175.8064} 
\email{omar.kchit@usmba.ac.ma, \,\, orcid: 0000.0002.0844.5034}
\begin{document}
	\begin{abstract}
		In this	paper, we deal with the problem of monogenity of number fields defined by monic irreducible trinomials $F(x)=x^{12}+ax^m+b\in \mathbb{Z}[x]$ with $1\leq m\leq11$. We give sufficient conditions on $a$, $b$, and $m$ so that the number field $K$ is not monogenic. In particular, for $m=1$ and for every rational prime $p$, we characterize when $p$ divides the index of $K$ and we provide a partial answer to the Problem $22$ of Narkiewicz \cite{Nar} for these number fields. Our results are illustrated by computational examples.
	\end{abstract}
	\maketitle
	\section{Introduction}
	Let $K=\mathbb{Q}(\alpha)$ be a number field of degree $n$ and $\mathbb{Z}_K$ its ring of integers. For any primitive element $\alpha\in \mathbb{Z}_K$ of $K$, it is well known that $\mathbb{Z}[\alpha]$ is a free $\mathbb{Z}$-module, which  occurs that $(\mathbb{Z}_K: \mathbb{Z}[\alpha])$ the index of  $\mathbb{Z}[\alpha]$ in $\mathbb{Z}_K$ is finite, and the formula linking $(\mathbb{Z}_K: \mathbb{Z}[\alpha])$, $\Delta(F)$, and $d_K$ is given by (\cite{G19}):
	\begin{equation}\tag{1.1}
		\Delta(\alpha)=\pm (\mathbb{Z}_K: \mathbb{Z}[\alpha])^2\times d_K,
	\end{equation}  
	where $d_K$ is the absolute discriminant of $K$ and $\Delta(\alpha)$ is the discriminant of the minimal polynomial of $\alpha$.  The greatest common divisor of the indices of all integral primitive elements of $K$ is called the  index of $K$, and  denoted by $i(K)$. Say $i(K)=\gcd \ \{ ( \mathbb{Z}_K; \mathbb{Z}[\alpha]) \, |\,\alpha \in \widehat{\mathbb{Z}}_K \}$. A rational prime $p$ dividing $i(K)$ is called a prime common index divisor of $K$. If $\mathbb{Z}_K$ has a power integral basis, i.e. a $\mathbb{Z}$-basis of the form $(1,\theta,\dots,\theta^{n-1})$, then $K$ is said to be monogenic. If $K$ is monogenic, then the index of $K$ is trivial, namely $i(K)=1$. Therefore a field having a prime common index divisor is not monogenic.  Problem $22$ of Narkiewicz \cite{Nar} asks for an explicit formula for the highest power of a given rational prime $p$ dividing $i(K)$, say $\nu_p(i(K))$. In 1871, Dedekind (\cite[§ 5, page 30]{R}) was the first who gave an example of a number field with non trivial index, he considered the cubic field $K$ generated by a root of $x^3-x^2-2x-8$ and showed that the rational prime $2$ splits completely in $K$. According to a theorem of Dedekind (\cite[Chapter I, Proposition 8.3]{Neu}), if we suppose that $K$ is monogenic, then we would be able to find a cubic polynomial generating $K$, that splits completely into distinct polynomials of degree $1$ in $\mathbb{F}_2[x]$. Since there are only two distinct polynomials of degree $1$ in $\mathbb{F}_2[x]$, this is impossible. Based on these ideas and using Kronecker's theory of algebraic number fields, Hensel  gave necessary and sufficient condition on the so-called "index divisors of $K$" for any rational prime $p$ to be a prime common index divisor \cite{He2}. The most important study of prime power decomposition of $i(K)$ has been done by Engstrom in 1930 (\cite{En}). He showed that $\nu_p(i(K))$ for number fields of degree $n\leq 7$ is determined by the factorization of  $p\mathbb{Z}_K$ into powers of prime ideals of $K$ for every positive prime $p\leq n$. In \cite{Sl}, \'Sliwa extended Engstom's results to number fields up to degree $12$ only for non-ramified primes $p$. These results were generalized by Nart \cite{N}, who developed a $p$-adic characterization of the index of a number field.  In \cite{Nak}, Nakahara  studied the  index of non-cyclic but abelian biquadratic number fields.  In \cite{GPP}  Ga\'al et al. characterized  the field indices of biquadratic number fields having Galois group $V_4$. In \cite{BM}, Bayad and Seddik studied the field index and determined the explicit prime ideal factorization of rational primes in a family of simplest quartic number fields. In \cite{DS}, for any  quartic number field $K$ defined by a trinomial $x^4+ax+b$, Davis and Spearman gave necessary and sufficient conditions on $a$ and $b$, which characterize when $p=2,3$ divides $i(K)$. Also in \cite{EG}, for any  quartic number field $K$ defined by a trinomial $x^4+ax^2+b$, El Fadil and Ga\'al gave necessary and sufficient conditions on $a$ and $b$, which characterize when a rational prime $p$ divides $i(K)$. In \cite{E5}, for any  quintic number field $K$ defined by a trinomial $x^5+ax^2+b$, El Fadil characterized when $p$ divides the index $i(K)$ for these number fields. In \cite{EK}, for any sextic number field defined by a trinomial $x^6+ax^5+b$, El Fadil and Kchit characterized $\nu_p(i(K))$ for $p=2,3,5$. In this paper, {for any number field} $K$ defined by a monic irreducible trinomial $F(x)=x^{12}+ax^m+b\in \mathbb{Z}[x]$, we give sufficient  conditions on $a$, $b$, and $m$ so that the number field $K$ is not monogenic. In particular, for $m=1$ and every rational prime $p$, we characterize when $p$ divides the index $i(K)$ and we provide a partial answer to the Problem $22$ of Narkiewicz \cite{Nar} for these number fields.
	\section{Main Results}
	Throughout this section, $K$ is a number field generated by a root $\alpha$ of a {monic irreducible} trinomial $F(x)=x^{12}+ax^m+b\in\mathbb{Z}[x]$ with $m=1,\dots,11$ is a natural integer. Without loss of generality, we assume that for every rational prime $p$,  $\nu_p(a)<12-m$ or $\nu_p(b)<12$.\\
	\smallskip
	
	The following theorem gives sufficient conditions on $a$ and $b$ so that $K$ is not monogenic, for every  $ m\in\{1,\dots,11\}$. 
	\begin{thm}\label{thm1}
		For every $m=1,\dots,11$, the following statements hold:
		\begin{enumerate}
			\item If $(a,b)\in \{(0,3),(0,7)\}\md{8}$, then $2$ divides $i(K)$.
			\item If $(a,b)\equiv(0,\pm1)\md{9}$, then $3$ divides $i(K)$.
		\end{enumerate}
		In particular, if one of the above conditions holds, then $K$ is not monogenic. 
	\end{thm}
	The following theorem gives sufficient conditions on $a$, $b$, and $m$ so that $K$ is not monogenic.
	\begin{thm}\label{thm2}
		If one of the following conditions holds:
		\begin{enumerate}
			\item $m\in\{3,4,8,9\}$ and $(a,b)\in\{(2,1),(6,5)\}\md{8}$,
			\item $m\in\{3,6,9\}$ and $(a,b)\equiv (4,7)\md{8}$,
		\end{enumerate}
		then $2$ divides $i(K)$. In particular, $K$ is not monogenic.
	\end{thm}
	The next theorem deals with a more general case of the family of number fields defined by the trinomials $x^n+ax^{2k}+b\in\mathbb{Z}$, where $n$ and $k$ are two positive integers such that $2k<n$. For every integer $a$, let $a_p=\cfrac{a}{p^{\nu_p(a)}}$.
	\begin{thm}\label{thm3}
		Let $K$ be a number field defined by a monic irreducible polynomial $F(x)=x^n+ax^{2k}+b\in\mathbb{Z}[x]$, where $n$ and $k$ are two positive integers $(2k<n)$. Let $p$ be an odd rational prime , if the following conditions hold:
		\begin{enumerate}
			\item $p\nmid (2k\times(n-2k))$,
			\item $p-1\mid (n-2k)$,
			\item $\nu_p(b)\geq1$ is even,
			\item $a\equiv -b_p\equiv - 1\md{p}$,
		\end{enumerate}
		then $p$ divides $i(K)$. In particular, $K$ is not monogenic.
	\end{thm}
	\begin{cor}\label{cor}
		Let $K$ be a number field defined by a monic irreducible trinomial $F(x)=x^{12}+ax^m+b$, then the following statements hold: 
		\begin{enumerate}
			\item If $m\in\{2,4,8,10\}$, $\nu_3(b)\geq 1$ is even, and $a\equiv -b_3\equiv -1\md{3}$, then $3$ divides $i(K)$.
			\item If $m\in\{4,8\}$, $\nu_5(b)\geq 1$ is even, and $a\equiv -b_5\equiv -1\md{5}$, then $5$ divides $i(K)$.
			\item If $m=6$, $\nu_7(b)\geq 1$ is even, and $a\equiv -b_7\equiv -1\md{7}$, then $7$ divides $i(K)$.
			\item If $m=2$, $\nu_{11}(b)\geq 1$ is even, and $a\equiv -b_{11}\equiv -1\md{11}$, then $11$ divides $i(K)$.
		\end{enumerate}
		In particular, if one of the above conditions holds, then $K$ is not monogenic.
	\end{cor}
	
	\smallskip In the next theorems, for $m=1$, and for every rational prime $p$, we characterize when $p$ divides the index $i(K)$.\\
	\smallskip
	{We start with the following Theorem which} characterizes when the ring $\mathbb{Z}[\alpha]$ is integrally closed. 
	\begin{thm}\label{mong}
		The ring $\mathbb{Z}[\alpha]$ is integrally closed if and only if the following conditions hold:
		\begin{enumerate}
			\item $b$ is square free.
			\item If $2$ divides $a$ and does not divide $b$, then $(a,b)\in\{(0,1),(2,3)\}\md{4}$.
			\item If $3$ divides $a$ and not divide $b$, then $(a,b)$ is in the set $$\{(0,2),(0,4),(0,5),(0,7),(3,-1),(3,1),(3,4),(3,7),(6,-1),(6,1),(6,4),(6,7)\}\md{9}.$$
			
			\item If $11$ divides $b$ and does not divide $a$, then $\nu_{11}(b-a^2+a^12)=1$.
			\item For every rational prime $p\not\in\{2,3,11\}$, if $\nu_p(ab)=0$, then $\nu_p(2^{24}\times 3^{12}b^{11}-11^{11}a^{12})\leq1$.
		\end{enumerate}
	\end{thm}
	\smallskip
	\begin{thm}\label{thmp2}
		Let $u=\nu_2(a+b+1)$ and $v=\nu_2(a+12)$. The rational prime $2$ divides the index $i(K)$ if and only if one of the following conditions holds:
		\begin{enumerate}
			\item $(a,b)\in \{(0,3),(0,7)\}\md{8}$. In particular, if $(a,b)\equiv (0,3)\md{8}$, then $\nu_2(i(K))=2$.
			\item $(a,b)\equiv (4,3)\md{8}$ and $u>2v-1$ or $u<2v-1$, $u$ is odd, and $\nu_2(b+as+s^{12})>2\nu_2(a+12s^{11})-1$ for some odd integer $s$. In this case, $\nu_2(i(K))=1$.
			\item $(a,b)\in\{(64,112),(192,112),(240,0),(240,128)\}\md{256}$.
			
			\item $(a,b)\in\{(256,832),(768,832),(0,64),(512,64)\}\md{1024}$. In particular, if $(a,b)\in\{(0,64),(512,64)\}\md{1024}$, then $\nu_2(i(K))=3$.
			\item $(a,b)\in\{(0,192),(256,192)\}\md{512}$. In this case, $\nu_2(i(K))=2$.
			\item $a\equiv 0\md{2^{11}}$ and $b\equiv 768+2^{10}+2^{11}\md{2^{12}}$.
		\end{enumerate}
		In particular, if one of the above conditions holds, then $K$ is not monogenic.
	\end{thm}
	\begin{thm}\label{thmp3}
		Let $u=\nu_3(b+a+1)$, $v=\nu_3(a+12)$, $\mu=\nu_3(b-a+1)$, and $\tau=\nu_3(a-12)$. The rational prime $3$ divides the index $i(K)$ if and only if one of the following conditions holds:
		\begin{enumerate}
			\item $(a,b)\equiv (0,\pm1)\md{9}$.
			
			\item $(a,b)\in\{(15,65),(42,38),(69,11)\}\md{81}$ and $u>2v-1$ or $u<2v-1$, $u$ is odd, and $(b+a+1)_3\equiv -1\md{3}$ or $u=2v-1$, $(b+a+1)_3\equiv 1\md{3}$, and $\nu_3(b+as+s^{12})>2\nu_3(a+12s^{11})-1$ or $\nu_3(b+as+s^{12})<2\nu_3(a+12s^{11})-1$, $\nu_3(b+as+s^{12})$ is odd, and $(b+as+s^{12})_3\equiv -1\md{3}$ for some integer $s$ such that $s\equiv 1\md{3}$.
			
			\item $(a,b)\in\{(6,47),(33,20),(60,74),(24,56),(51,29),(78,2)\}\md{81}$ and $\nu_3(b+as+s^{12})>2\nu_3(a+12s^{11})$ or $\nu_3(b+as+s^{12})<2\nu_3(a+12s^{11})$, $\nu_3(b+as+s^{12})$ is odd, and $(b+as+s^{12})_3\equiv -1\md{3}$ for some integer $s$ such that $s\equiv 1\md{3}$.
			
			\item $(a,b)\in\{(12,11),(39,38),(66,65)\}\md{81}$, and $\mu>2\tau-1$ or $\mu<2\tau-1$, $\mu$ is odd, and $(b-a+1)_3\equiv -1\md{3}$ or $\mu=2\tau-1$, $(b-a+1)_3\equiv 1\md{3}$, and $\nu_3(b+as+s^{12})>2\nu_3(a+12s^{11})-1$ or $\nu_3(b+as+s^{12})<2\nu_3(a+12s^{11})-1$, $\nu_3(b+as+s^{12})$ is odd, and $(b+as+s^{12})_3\equiv -1\md{3}$ for some integer $s$ such that $s\equiv -1\md{3}$.
			
			\item $(a,b)\in\{(3,2),(30,29),(57,56),(21,74),(48,20),(75,47)\}\md{81}$ and $\nu_3(b+as+s^{12})>2\nu_3(a+12s^{11})-1$ or $\nu_3(b+as+s^{12})<2\nu_3(a+12s^{11})-1$, $\nu_3(b+as+s^{12})$ is odd, and $(b+as+s^{12})_3\equiv -1\md{3}$ for some integer $s$ such that $s\equiv -1\md{3}$.
		\end{enumerate}
		In particular, if one of the above conditions holds, then $K$ is not monogenic.\\
		Furthermore, except the case $(1)$, we have $\nu_3(i(K))=1$.
	\end{thm}
	\begin{thm}\label{pro1}
		For every rational prime $p\geq 5$ and every $(a,b)\in\mathbb{Z}^2$, $p$ does not divide the index $i(K)$.
	\end{thm}
	
	\section{Preliminaries}
	
	{Since our proofs are based on prime ideal factorization and} Newton polygon techniques is a standard method which is rather technical but very efficient to apply. We briefly describe the use of these techniques, which makes our proofs understandable. For more details, we refer to \cite{EMN} and \cite{GMN}.\\
	Let $K=\mathbb{Q}(\alpha)$ be a number field generated by a  root $\alpha$ of a monic irreducible polynomial $F(x)\in\mathbb{Z}[x]$. We shall use Dedekind’s theorem \cite[25, Chapter I, Proposition 8.3]{Neu} relating the prime
	ideal factorization of $p\mathbb{Z}_K$ and the factorization of $F(x)$ modulo $p$ $($for rational primes $p$ not
	dividing $(\mathbb{Z}_K:\mathbb{Z}[\alpha]))$. Also, we shall need Dedekind’s criterion \cite[5, Theorem 6.1.4]{Co}
	on the divisibility of $(\mathbb{Z}_K:\mathbb{Z}[\alpha])$ by rational primes $p$.\\
	For any rational prime $p$, let $\nu_p$ be the $p$-adic valuation of $\mathbb{Q}$, $\mathbb{Q}_p$ its $p$-adic completion, and $\mathbb{Z}_p$ the ring of $p$-adic integers. We denote by $\nu_p$ be the Gauss's extension of $\nu_p$ to $\mathbb{Q}_p(x)$ and is defined on $\mathbb{Q}_p[x]$ by $\nu_p(\sum_{i=0}^na_ix^i)=\min_i\{\nu_p(a_i)\},a_i\in\mathbb{Q}_p$. Also, for nonzero polynomials, $P,Q\in\mathbb{Q}_p[x]$, we extend this valuation to $\nu_p(P/Q)=\nu_p(P)-\nu_p(Q)$. Let $\phi\in\mathbb{Z}_p[x]$ be a monic lift to an irreducible factor of $F(x)$ modulo $p$. Upon the Euclidean division by successive powers of $\phi$, there is a unique {$\phi$-expansion} of $F(x)$; that is $F(x)=a_0(x)+a_1(x)\phi(x)+\cdots+a_l(x)\phi(x)^l$, where $a_i(x)\in\mathbb{Z}_p[x]$ and deg$(a_i)<$deg$(\phi)$. For every $i=0,\dots,l$, let  $u_i=\nu_p(a_i(x))$. The {$\phi$-Newton polygon} of $F(x)$ with respect to $p$, is the lower boundary convex envelope of the set of points $\{(i,u_i),~ a_i(x)\neq 0\}$ in the Euclidean plane, which we denote by $N_{\phi}(F)$. It is the process of joining the obtained edges $S_1,\dots,S_t$ ordered by increasing slopes, which can be expressed as $N_{\phi}(F)= S_1+\cdots+S_t$. For every
	side $S_i$ of $N_{\phi}(F)$, its length $l(S_i)$ is the length of its projection to the $x$-axis and its
	height $h(S_i)$ is the length of its projection to the $y$-axis. We call $d(S_i)=$gcd$(l(S_i),h(S_i))$ the degree of $S_i$. The polygon determined by the sides of the $\phi$-Newton polygon with negative slopes is called the {principal $\phi$-Newton polygon} of $F(x)$, and it is denoted by $N_{\phi}^+(F)$. As defined in \cite[Def. 1.3]{EMN}, the $\phi$-{index} of $F(x)$, denoted  $ind_{\phi}(F)$, is  deg$(\phi)$ multiplied by the number of points with natural integer coordinates that lie below or on the polygon $N_{\phi}^+(F)$, strictly above the horizontal axis,{ and strictly beyond the vertical axis}. Let $\mathbb{F}_{\phi}$ be the field $\mathbb{F}_p[x]/(\ol{\phi})$, then to every side $S$ of $N_{\phi}^+(F)$ with initial point $(i,u_i)$, and every $i=0,\ldots,l$, let the residue coefficient $c_i\in\mathbb{F}_{\phi}$ defined as follows: 
	$$c_{i}=
	\left
	\{\begin{array}{ll} 0,& \mbox{ if } (s+i,{\it u_{s+i}}) \mbox{ lies strictly
			above } S,\\
		\left(\dfrac{a_{s+i}(x)}{p^{{\it u_{s+i}}}}\right)
		\,\,
		\mod{(p,\phi(x))},&\mbox{ if }(s+i,{\it u_{s+i}}) \mbox{ lies on }S.
	\end{array}
	\right.$$
	where $(p,\phi(x))$ is the maximal ideal of $\mathbb{Z}[x]$ generated by $p$ and $\phi$. Let $\lambda=-h/e$ be the slope of $S$, where  $h$ and $e$ are two positive coprime integers and $l=l(S)$. Then  $d=l/e$ is the degree of $S$. Since  the points  with integer coordinates lying{ on} $S$ are exactly ${(s,u_s),(s+e,u_{s}-h),\dots, (s+de,u_{s}-dh)}$. Thus if $i$ is not a multiple of $e$, then 
	$(s+i, u_{s+i})$ does not lie on $S$, and so $c_i=0$. Let
	${R_{\lambda}(F)(y)}=t_dy^d+t_{d-1}y^{d-1}+\cdots+t_{1}y+t_{0}\in\mathbb{F}_{\phi}[y]$, called  
	the {residual polynomial} of $F(x)$ associated to the side $S$, where for every $i=0,\dots,d$,  $t_i=c_{ie}$. If ${R_{\lambda}(F)(y)}$ is square free for each side of the polygon $N_{\phi}^+(F)$, then we say that $F$ is $\phi$-regular. Let $\ol{F(x)}=\prod_{i=1}^{r}\ol{\phi_i}^{l_i}$ be the factorization of $F(x)$ into powers of monic irreducible coprime polynomials over $\mathbb{F}_p$, we say that the polynomial $F(x)$ is $p$-{regular} if $F(x)$ is a $\phi_i$-regular polynomial with respect to $p$ for every $i=1,\dots,r$. Let  $N_{\phi_i}^+=S_{i1}+\cdots+S_{ir_i}$ be the $\phi_i$-principal Newton polygon of $F(x)$ with respect to $p$. For every $j=1,\dots,r_i$, let let $R_{\lambda_{ij}}(y)=\prod_{s=1}^{s_{ij}}\psi_{ijs}^{a_{ijs}}(y)$ be the factorization of $R_{\lambda_{ij}}(y)$ in $\mathbb{F}_{\phi_i}[y]$. Then we have the following  theorem of index of Ore:
	\begin{thm}\label{thm4}$($\cite[Theorem 1.7 and Theorem 1.9]{EMN}$)$\\
		Under the above hypothesis, we have the following:
		\begin{enumerate}
			\item 
			$$\nu_p((\mathbb{Z}_K:\mathbb{Z}[\alpha]))\geq\sum_{i=1}^{r}ind_{\phi_i}(F).$$  
			The equality holds if $F(x) \text{ is }p$-regular.
			\item 
			If  $F(x) \text{ is }p$-regular, then
			$$p\mathbb{Z}_K=\prod_{i=1}^r\prod_{j=1}^{t_i}\prod_{s=1}^{s_{ij}}\mathfrak{p}_{ijs}^{e_{ij}}.$$
			{where $e_{ij}$ is the smallest positive integer satisfying $e_{ij}\la_{ij}\in \Z$ and\\
				$f_{ijs}=\deg(\phi_i)\times \deg(\psi_{ijs})$ is the residue degree of $\mathfrak{p}_{ijs}$ over $p$ for every $(i,j,s)$}.
		\end{enumerate}
	\end{thm}
	\begin{cor}\label{cor1}
		Under the hypothesis  above {$($Theorem $\ref{thm4})$}, if for every $i=1,\dots,r,\,l_i=1\text{ or }N_{\phi}^+(F)=S_i$ has a single side of height $1$, then $\nu_p((\mathbb{Z}_K:\mathbb{Z}[\alpha]))=0$.
	\end{cor}
	If a factor of $F(x)$ provided by Hensel’s lemma and refined by  Newton polygon (in the context of Ore {program}) is  not irreducible over $\mathbb{Q}_p$, then in order to complete the factorization of $F(x)$ in $\mathbb{Q}_p[x]$, Guardia, Montes, and Nart introduced the notion of {higher order Newton polygon} {\cite{GMN}}. They showed,  thanks to a theorem of index \cite[Theorem 4.18]{GMN}, that  after a finite number of iterations the process provides all monic irreducible factors of $F(x)$ {and all prime ideals of $\mathbb{Z}_K$ lying above a rational prime $p$}. We recall some fundamental techniques of Newton polygon of high order. For more details, we refer to \cite{GMN} and \cite{GN}. A {type of order} $r - 1$ is  a data ${\bf{t}}  = (g_1(x), -\la_1, g_2(x), -\la_2,\dots, g_{r-1}(x), -\la_{r-1},\psi_{r-1}(x))$,
	where every $g_{i}(x)$ is a monic polynomial in $\mathbb{Z}_p[x]$, $\la_i\in \mathbb{Q}^+$ and
	$\psi_{r-1}(y)$ is a polynomial over a  finite field of $\displaystyle p^{H}$ elements,  with  $H=\displaystyle\prod_{i=0}^{r-2}f_i$ and $f_i=\mbox{deg}(\psi_{i}(x))$, satisfying the following recursive properties:
	\begin{enumerate}
		\item[(0)] {$\mathbb{F}_0$ is the finite field of $p$ elements}.
		\item[(1)]
		$g_{1}(x)$ is irreducible modulo $p$, {$\psi_0(y) \in \mathbb{F}_0[y]$} {is} the polynomial obtained
		by reducing $g_{1}(x)$ modulo $p$, and  $\mathbb{F}_1 = \mathbb{F}_0[y]/{(\psi_0(y))}$.
		\item[(2)]
		For every $i=1,\dots,r-1$, the Newton polygon of $i^{th}$ order, $N_i(g_{i+1}(x))$,  has a single  side of slope $-\la_i$.
		\item[(3)] 
		For every $i=1,\dots,r-1$, the residual polynomial of $i^{th}$ order, $R_i(g_{i+1})(y)=\psi_i(y)\in \mathbb{F}_i[y]$ is a monic irreducible polynomial in $\F_i[y]$, and  $\F_{i+1}= \F_i[y]/(\psi_i(y))$.
		\item[(4)]
		For every $i=1,\dots,r-1$, $g_{i+1}(x)$ has minimal degree among all monic polynomials
		in $\mathbb{Z}_p[x]$ satisfying $(2)$ and $(3)$.
		\item[(5)] 
		$\psi_{r-1}(y) \in\F_{r-1}[y]$ is a monic irreducible polynomial,  {$\psi_{r-1}(y)\neq y$}, and $\F_{r}= \F_{r-1}[y]/(\psi_{r-1}(y))$.
	\end{enumerate}
	Thus $\F_0\subset \F_1\subset \dots\subset \F_r$ is a tower of finite fields, here  the field
	$\F_i$ should not be confused with the finite field of $i$ elements. Let $\om_0=[\nu_p,x,0]$ be the Gauss's extension of $\nu_p$ to $\mathbb{Q}_p(x)$. {Since $R_i(g_{i+1})(y)$ is an irreducible polynomial in $\F_i[y]$, and by the theorem of the product in order $i$, the polynomial $g_i(x)$ is irreducible in $\mathbb{Z}_p[x]$. Then according} to MacLane notations  {and definitions} (\cite{Mc}),  $g_{i+1}$ is a key polynomial of $\om_i$, and it induces  a valuation  on $\Q_p(x)$,  denoted by  $\om_{i+1}=e_i[\om_{i},g_{i+1},\la_{i+1}]$, where $\la_i=h_i/e_i$,   $e_i$ and  $h_i$ are  positive coprime integers. The valuation   $\om_{i+1}$ is called the augmented valuation of $\nu_p$  with respect to $\ph$ and $\la$ is  defined over   $\Q_p[x]$  as follows:   
	
	$$\om_{i+1}(F(x))=\mbox{min}\left\{e_{i+1}\om_i(a_j^{i+1}(x))+jh_{i+1},\,j=0,\dots,n_{i+1}\right\},$$ 
	where $F(x)=\displaystyle\sum_{j=0}^{n_{i+1}}a_j^{i+1}(x)g_{i+1}^j(x)$ is the $g_{i+1}$-expansion of $F(x)$. According to the terminology  in \cite{GMN}, the valuation  $\om_r$ is called the  $r^{th}$-order valuation associated to the data ${\bf{t}}$. For every order $r\ge 1$, the $g_r$-Newton polygon of {$F(x)$} with respect  to  the valuation $\om_r$ is the lower boundary convex envelope of  the set of  points $\{(i,\mu_i), i=0,\dots, n_r\}$ in the  {Euclidean} plane, where $\mu_i=\om_r(a^r_i(x)g_r^i(x)))$. {Theorem of the product, theorem of the polygon, and theorem of the residual polynomial in high order Newton polygon are the relevant theorems from Montes-Guardia-Nart's work (see \cite[Theorems 2.26, 3.1, 3.7]{GMN})}.\\
	\smallskip
	For the proof of our main results, we need the following lemma, which characterizes the prime common index divisors of $K$.
	\begin{lem}\label{index}$($\cite{En}$)$\\
		Let $p$ be a rational prime and $K$ be a number field. For every positive integer $f$, let $\mathcal{P}_f$ be the number of distinct prime ideals of $\mathbb{Z}_K$ lying above $p$ with residue degree $f$ and $\mathcal{N}_f$ the number of monic irreducible polynomials of $\mathbb{F}_p[x]$ of degree $f$.  Then $p$ {divides the index $i(K)$} if and only if $\mathcal{P}_f>\mathcal{N}_f$ for some positive integer $f$.
	\end{lem}
		
		\section{Proofs of Main Results}
		Throughout this section, $p$ is a rational prime, $\mathfrak{p}^e$ is a prime ideal of $\mathbb{Z}_K$ lying above $p$, $f$ its residue degree, and $e$ its ramification index. For every integer $a$, let $a_p=\cfrac{a}{p^{\nu_p(a)}}$.\\

		\textit{\textbf{ Proof of Theorem \ref{thm1}}}.\\
		In every case, let us show that $i(K)>1$, and so $K$ is not
		monogenic.
		\begin{enumerate}
			\item For $(a,b)\in\{(0,3),(0,7)\}\md{8}$, we have $F(x)\equiv (x-1)^4(x^2+x+1)^4\md{2}$ for every $ m=1,\dots,11$. Let $\phi_1=x-1$ and $\phi_2=x^2+x+1$, then by Taylor expansion we get $$
			\begin{array}{lll}
				F(x)&=&\cdots+\left(495+\cfrac{am(m-1)(m-2)(m-3)}{24}\right)\phi_1^4+\left(220+\cfrac{a(m-1)(m-2)}{6}\right)\phi_1^3\\
				&+&\left(66+\cfrac{am(m-1)}{2}\right)\phi_1^2+(12+am)\phi_1+a+b+1.
			\end{array}
			$$
			
			The $\phi_2$-expansion of $F(x)$ are given by the following equations corresponding to the values of $m$.
			$$
			\begin{array}{llll}
				m=1:&F(x)&=&\cdots+(-25-5x)\phi_2^4+(18+24x)\phi_2^3-18x\phi_2^2+(-4+4x)\phi_2+ax+b+1,\\
				m=2:&F(x)&=&\cdots+(-25-5x)\phi_2^4+(18+24x)\phi_2^3-18x\phi_2^2+((a+4)x-a-4)\phi_2\\
				&&&-ax-a+b+1,\\
				m=3:&F(x)&=&\cdots+(-25-5x)\phi_2^4+(18+24x)\phi_2^3-18x\phi_2^2+((a+4)x-a-4)\phi_2\\
				&&&+a+b+1,\\
				m=4:&F(x)&=&\cdots+(-25-5x)\phi_2^4+(18+24x)\phi_2^3+(-18x+a)\phi_2^2\\
				&&&+((-2a+4)x-a-4)\phi_2+ax+b+1,\\				
				m=5:&F(x)&=&\cdots+(-25-5x)\phi_2^4+(18+24x)\phi_2^3+((a-18)x-2a)\phi_2^2\\
				&&&+((a+4)x+3a-4)\phi_2-ax-a+b+1,\\
				m=6:&F(x)&=&\cdots+(-25-5x)\phi_2^4+(18+24x+a)\phi_2^3+((-3a-18)x)\phi_2^2\\
				&&&+((2a+4)x-2a-4)\phi_2+a+b+1,\\
				m=7:&F(x)&=&\cdots+((a-5)x-4a-25)\phi_2^4+((a+24)x-3a+18)\phi_2^3\\
				&&&+((-18+3a)x+5a)\phi_2^2+((-4a+4)x-2a-4)\phi_2+ax+b+1,\\
				m=8:&F(x)&=&\cdots+(-25+a-5x)\phi_2^4+((-4a+24)x+18+2a)\phi_2^3\\
				&&&+((2a-18)x-7a)\phi_2^2+((2a+4)x+5a-4)\phi_2-ax-a+b+1,\\
				m=9:&F(x)&=&\cdots+((a-5)x-4a-25)\phi_2^4+((6a+24)x+6a+18)\phi_2^3\\
				&&&+((-9a-18)x)\phi_2^2+((3a+4)x-3a-4)\phi_2+a+b+1,\\
				m=10:&F(x)&=&\cdots+((-5a-5)x+5a-25)\phi_2^4+(18+24x-15a)\phi_2^3\\
				&&&+((-18+9a)x+12a)\phi_2^2+((-6a+4)x-3a-4)\phi_2+ax+b+1,\\
				m=11:&F(x)&=&\cdots+((10a-5)x+5a-25)\phi_2^4+((-15a+24)x+18+9a)\phi_2^3\\
				&&&+((-18+3a)x-15a)\phi_2^2+((3a+4)x+7a-4)\phi_2-ax-a+b+1.
			\end{array}
			$$
			
			\begin{enumerate}
				\item[(i)] If $(a,b)\equiv (0,3)\md{8}$, then $N_{\phi_2}^+(F)=S_2$ has a single side joining $(0,2)$ and $(2,0)$ with $R_{\lambda_2}(F)(y)=(x+1)y^2+xy+1=(y+1)(y+x)\in\mathbb{F}_{\phi_2}[y]$. By Theorem $\ref{thm4}$, $\phi_2$ provides two prime ideals of $\mathbb{Z}_{K}$ lying above $2$ with residue degree $2$ each. Since there is only one monic irreducible polynomial of degree $2$ in $\mathbb{F}_2[x]$ namely $x^2+x+1$, by Lemma $\ref{index}$, $2$ divides $i(K)$ and so $K$ is not monogenic.
				
				\item[(ii)] If $(a,b)\equiv (0,7)\md{8}$, then for $(a,b)\in\{(0,-1),(8,7)\}\md{16}$ we have $w=\nu_2(a+b+1)\geq 4$, $\nu_2(12+am)=2$, and $\nu_2\left(66+\cfrac{am(m-1)}{2}\right)=1$ for every $m=1,\dots,11$. Thus $N_{\phi_1}^+(F)=S_1+S_2+S_3$ has three sides joining $(0,w)$, $(1,2)$, $(2,1)$, and $(4,0)$. By Theorem $\ref{thm4}$, $\phi_1$ provides three prime ideals of $\mathbb{Z}_{K}$ lying above $2$ with residue degree $1$ each. Since there is just two monic irreducible polynomials of degree $1$ in $\mathbb{F}_2[x]$, by Lemma $\ref{index}$, $2$ divides $i(K)$ and so $K$ is not monogenic. For $(a,b)\in\{(0,7),(8,-1)\}\md{16}$ we have $\nu_2(a+b+1)=3$, $\nu_2(12+am)=2$, and $\nu_2\left(66+\cfrac{am(m-1)}{2}\right)=1$. Thus $N_{\phi_1}^+(F)=S_{11}+S_{12}$ has two sides joining $(0,3)$, $(2,1)$, and $(4,0)$ with $d(S_{12})=1$ and $R_{\lambda_{11}}(F)(y)=y^2+y+1$, which is irreducible over $\mathbb{F}_{\phi_1}$. Hence, $\phi_1$ provides two prime ideals of $\mathbb{Z}_{K}$ lying above $2$ with residue degrees $1$ and $2$ respectively. On the other hand, $N_{\phi_2}^+(F)=S_{21}+S_{22}$ has two sides joining $(0,3)$, $(2,1)$, and $(4,1)$ with $d(S_{22})=1$. Thus $\phi_2$ provides at least one prime ideal of $\mathbb{Z}_{K}$ lying above $2$ with residue degree $2$. We conclude that there are at least two prime ideals of $\mathbb{Z}_{K}$ lying above $2$ with residue degree $2$ each. By Lemma $\ref{index}$, $2$ divides $i(K)$ and so $K$ is not monogenic.
			\end{enumerate}	
			\item For $(a,b)\equiv (0,\pm1)\md{9}$, we have the following cases:
			\begin{enumerate}
				\item[(i)] If $b\equiv -1\md{9}$, then $F(x)\equiv (x-1)^3(x+1)^3(x^2+1)^3\md{3}$. Let $\phi_1=x-1$ and $\phi_2=x+1$, then by Taylor expansion we get 
				$$
				\begin{array}{lll}
					F(x)&=&\cdots+\left(495+\cfrac{am(m-1)(m-2)(m-3)}{24}\right)\phi_1^4+\left(220+\cfrac{a(m-1)(m-2)}{6}\right)\phi_1^3\\
					&+&\left(66+\cfrac{am(m-1)}{2}\right)\phi_1^2+(12+am)\phi_1+a+b+1,\\
					&=&\cdots+\left(-220+\cfrac{(-1)^{m-3}am(m-1)(m-2)}{6}\right)\phi_2^3+\left(66+\cfrac{(-1)^{m-2}am(m-1)}{2}\right)\phi_2^2\\
					&+&(-12+(-1)^{m-1}am)\phi_2+(-1)^ma+b+1.
				\end{array}$$
				We have $\nu_3((-1)^ma+b+1)\geq2$, $\nu_3(a+b+1)\geq2$, and $\nu_3(-12+(-1)^{m-1}am)=\nu_3(12+am)=1$ for every $m=1,\dots,11$. Thus $N_{\phi_i}^+(F)=S_{i1}+S_{i2}$ has two sides of degree $1$ each for every $i=1,2$. Hence, $\phi_i$ provides two prime ideals of $\mathbb{Z}_{K}$ lying above $3$ with residue degree $1$ each. Applying this for $i=1,2$, we conclude that there are four prime ideals of $\mathbb{Z}_{K}$ lying above $3$ with residue degree $1$ each. Since there is just three monic irreducible polynomials of degree $1$ in $\mathbb{F}_3[x]$, by Lemma $\ref{index}$, $3$ divides $i(K)$ and so ${K}$ is not monogenic.
				
				\item[(ii)] If $b\equiv 1\md{9}$, then $F(x)\equiv (x^2+x-1)^3(x^2-x-1)^3\md{3}$. Let $\phi_{\pm}=x^2\pm x-1$ $(\phi_+=x^2+x-1$ and $\phi_-=x^2-x-1)$, then the $\phi_{\pm}$-expansion of $F(x)$ are given by the following equations corresponding to the values of $m$.
				$$
				\begin{array}{llll}
					m=1:&F(x)&=&\cdots+(338\mp256x)\phi_{\pm}^3+(468\mp474x)\phi_{\pm}^2+(324\mp420x)\phi_{\pm}\\
					&&&+(a\mp144)x+b+89,\\
					
					m=2:&F(x)&=&\cdots+(338\mp256x)\phi_{\pm}^3+(468\mp474x)\phi_{\pm}^2+(\mp420x+a+324)\phi_{\pm}\\
					&&&\mp(144a+144)x+a+b+89,\\
					
					m=3:&F(x)&=&\cdots+(338\mp256x)\phi_{\pm}+(468\mp474x)\phi_{\pm}^2+((a\mp420)x+324\mp a)\phi_{\pm}\\
					&&&+(2a\mp144)x\mp a+b+89,\\
					m=4:&F(x)&=&\cdots+(338\mp256x)\phi_{\pm}^3+(\mp474x+a+468)\phi_{\pm}^2\\
					&&&+(\mp(420+2a)x+3a+324)\phi_{\pm}\mp(144+3a)x+2a+b+89,\\
				\end{array}
				$$

				$$
				\begin{array}{llll}
					m=5:&F(x)&=&\cdots+(338\mp256x)\phi_{\pm}^3+((a\mp474)x+468\mp2a)\phi_{\pm}^2\\
					&&&+((5a\mp420)x+324\mp5a)\phi_{\pm}+(5a\mp144)x\mp3a+b+89,\\
					m=6:&F(x)&=&\cdots+(\mp256x+338+a)\phi_{\pm}^3+(\mp(474+3a)x+6a+468)\phi_{\pm}^2\\
					&&&+(\mp(420+10a)x+10a+324)\phi_{\pm}
					\mp(144+8a)x+5a+b+89,\\
					m=7:&F(x)&=&\cdots+((a\mp256)x+338\mp3a)\phi_{\pm}^3+((9a\mp474)x+468\mp13a)\phi_{\pm}^2\\
					&&&+((20a\mp420)x+324\mp18a)\phi_{\pm}+(13a\mp144)x\mp8a+b+89,\\
					m=8:&F(x)&=&\cdots+(\mp(256+4a)x+10a+338)\phi_{\pm}^3+(\mp(474+22a)x+29a+468)\phi_{\pm}^2\\
					&&&+(\mp(420+38a)x+33a+324)\phi_{\pm}\mp(144+21a)x+13a+b+89,\\
					m=9:&F(x)&=&\cdots+((14a\mp256)x+338\mp26a)\phi_{\pm}^3+((51a\mp474)x+468\mp60a)\phi_{\pm}^2\\
					&&&+((71a\mp420)x+324\mp59a)\phi_{\pm}+(34a\mp144)x\mp21a+b+89,\\
					m=10:&F(x)&=&\cdots+(\mp(256+40a)x+65a+338)\phi_{\pm}^3+(\mp(474+111a)x+122a+468)\phi_{\pm}^2\\
					&&&+(\mp(420+130a)x+105a+324)\phi_{\pm}+\mp(144+55a)x+34a+b+89,\\
					m=11:&F(x)&=&\cdots+((105a\mp256)x+338\mp151a)\phi_{\pm}^3+((233a\mp474)x+468\mp241a)\phi_{\pm}^2\\
					&&&+((235a\mp420)x+324\mp185a)\phi_{\pm}+(89a\mp144)x\mp55a+b+89.
					
				\end{array}$$

				For every $m=1,\dots,11$, we have $N_{\phi_{\pm}}^+(F)=S_1+S_2$ has two sides joining $(0,u_i)$, $(1,1)$, and $(3,0)$ with $u_i\geq2$ for every $i=1,\dots,11$. Thus the degree of each side $1$. Hence $\phi_{\pm}$ provides two prime ideals of $\mathbb{Z}_{K}$ lying above $3$ with residue degree $2$ each. Applying that for $\phi_-$ and $\phi_+$, we get $3\mathbb{Z}_{K}=(2,2)(2,2)(2,1)(2,1)$. Since there is just three monic irreducible polynomials of degree $2$ in $\mathbb{F}_3[x]$ namely $x^2+1$, $x^2-x-1$, and $x^2+x-1$, by Lemma $\ref{index}$, $3$ divides $i(K)$ and so ${K}$ is not monogenic.  
			\end{enumerate}
		\end{enumerate}
		\begin{flushright}
			$\square$
		\end{flushright}
		
		\textit{\textbf{Proof of Theorem \ref{thm2}}}.\\
		\begin{enumerate}
			\item If $(a,b)\in\{(2,1),(6,5)\}\md{8}$, then $F(x)\equiv \phi_1^4\cdot\phi_2^4\md{2}$ with $\phi_1=x-1$ and $\phi_2=x^2+x+1$. The expansions of $F(x)$ are given above. For $m=4,8$, we have $\nu_2(a+b+1)=2$, $\nu_2(12+am)\geq3$, and  $\nu_2\left(66+\cfrac{am(m-1)}{2}\right)=1$. Thus $N_{\phi_1}^+(F)=S_1$ has a single side joining $(0,2)$ and $(4,0)$ with $R_{\lambda_1}(F)(y)=y^2+y+1$, which is irreducible over $\mathbb{F}_{\phi_1}$. Hence, $\phi_1$ provides a unique prime ideal of $\mathbb{Z}_{K}$ lying above $2$ with residue degree $2$. On the other hand, $N_{\phi_2}^+(F)=S_2$ has a single side of height $1$. Thus $2\mathbb{Z}_{K}=\mathfrak{p}_1^2\mathfrak{p}_2^4$ with residue degree $2$ each. Hence $2$ divides $i(K)$ and so $K$ is not monogenic. For $m=3,9$, we have $N_{\phi_2}(F)=S_{21}+S_{22}$ has two sides joining $(0,2)$, $(1,1)$, and $(0,4)$. Thus $\phi_2$ provides two prime ideals of $\mathbb{Z}_{K}$ lying above $2$ with residue degree $2$ each. Hence, $2$ divides $i(K)$ and so $K$ is not monogenic. 
			
			\item If $(a,b)\equiv(4,7)\md{8}$, then $F(x)\equiv \phi_1^4\cdot \phi_2^4\md{2}$ with $\phi_1=x-1$ and $\phi_2=x^2+x+1$. For $m\in\{3,6,9\}$, we have $N_{\phi_2}^+(F)=S_2$ has a single side joining $(0,2)$ and $(2,0)$ with $R_{\lambda_2}(F)(y)=(x+1)y^2+xy+1=(y+1)(y+x)\in\mathbb{F}_{\phi_2}[y]$. Thus $\phi_2$ provides two prime ideals of $\mathbb{Z}_{K}$ lying above $2$ with residue degree $2$ each. Hence, $2$ divides $i(K)$ and so ${K}$ is not monogenic.
		\end{enumerate}
		\smallskip
		
		\textit{\textbf{Proof of Theorem \ref{thm3}}}.\\
		If $\nu_p(b)\geq 1$ and $a\equiv -1\md{p}$, then $F(x)\equiv x^{2k}(x^{n-2k}-1)\md{p}$. Since $n-2k\not\equiv 0\md{p}$, then the polynomial $x^{n-2k}-1$ is separable over $\mathbb{F}_p$. Also, we have $p-1$ divides $n-2k$, then $(x^{n-2k}-1)\equiv (x^{p-1}-1)U(x)\md{p}$ such that $GCD(\overline{x^{p-1}-1},\overline{U(x)})=1$. Thus $F(x)\equiv x^{2k}(x^{p-1}-1)U(x)\equiv x^{2k}U(x)\prod_{j=1}^{p-1}(x-j)\md{p}$. Clearly, $(x-j)$ provides a unique prime ideal of $\mathbb{Z}_K$ lying above $p$ with residue degree $1$ each, for every $j=1,\dots,p-1$. On the other hand, let $\phi=x$, then $N_{\phi}^+(F)=S$ has a single joining $(0,\nu_p(b))$ and $(2k,0)$ with $R_{\lambda}(F)(y)=y^d-1$ where $d=gcd(2k,\nu_2(b))$. Since $p\nmid 2k$, then $R_{\lambda}(F)(y)$ is separable over $\mathbb{F}_{\phi}$. Finally, if $\nu_p(b)$ is even, then $2$ divides $d$, and so $-y^2+1=-(y-1)(y+1)$ divides $R_{\lambda}(F)(y)$. Thus $\phi$ provides at least two prime ideals of $\mathbb{Z}_K$ lying above $p$ with residue degree $1$ each. We conclude that there are at least $p+1$ prime ideals of $\mathbb{Z}_K$ lying above $p$ with residue degree $1$ each. Since there is just $p$ monic irreducible polynomial of degree $1$ in $\mathbb{F}_p[x]$, by Lemma $\ref{index}$, $p$ divides $i(K)$ and so $K$ is not monogenic.
		\begin{flushright}
			$\square$
		\end{flushright}
		\smallskip
		
		\textit{\textbf{Proof of Theorem \ref{mong}}}.\\
		Let $F(x)=x^{12}+ax+b$ and $p$ a rational prime candidate to divide $(\mathbb{Z}_K:\mathbb{Z}[\alpha])$. Since $\Delta(F)=2^{24}\times 3^{12}b^{11}-11^{11}a^{12}$ is the discriminant of $F(x)$ and thanks to the index formula $(1.1)$, we have that $p^2$ divides $\Delta(F)$; that is $p^2$ divides $2^{24}\times 3^{12}b^{11}-11^{11}a^{12}$. Hence, if $p=2,3$, then $p$ divides $a$. If $p=11$, then $p$ divides $b$. If $p\not\in\{2,3,11\}$, then $p$ divides $a$ and $b$ or $p$ does not divide $ab$ and $p^2$ divides $2^{24}\times 3^{12}b^{11}-11^{11}a^{12}$.
		\begin{enumerate}
			\item For $p=2$, assume that $2$ divides $a$.
			\begin{enumerate}
				\item[(a)] If $2$ divides $b$, then $2$ does not divide $(\mathbb{Z}_K:\mathbb{Z}[\alpha])$ if and only if $\nu_2(b)=1$.
				
				\item[(b)] If $2$ does not divide $b$, then $F(x)\equiv (x-1)^4(x^2+x+1)^4\md{2}$. Let $\phi_1=x-1$ and $\phi_2=x^2+x+1$, then 
				$$\begin{array}{lll}
					F(x)&=&\cdots+(a+12)\phi_1+a+b+1,\\
					&=&\cdots+(-4+4x)\phi_2+ax+b+1.
				\end{array}$$
				We conclude that $2$ does not divide $(\mathbb{Z}_K:\mathbb{Z}[\alpha])$ if and only if $\nu_2(a+b+1)=\nu_2(ax+b+1)=1$; $(a,b)\in\{(0,1),(2,3)\}\md{4}$.
			\end{enumerate}
			\item For $p=3$, assume that $3$ divides $a$.
			\begin{enumerate}
				\item[(a)] If $3$ divides $b$, then $3$ does not divide $(\mathbb{Z}_K:\mathbb{Z}[\alpha])$ if and only if $\nu_3(b)=1$.
				
				\item[(b)] If $b\equiv -1\md{3}$, then $F(x)\equiv (x-1)^3(x+1)^3(x^2+1)^3\md{3}$. Let $\phi_1=x-1$, $\phi_2=x+1$, and $\phi_3=x^2+1$, then 
				
				$$\begin{array}{lll}
					F(x)&=&\cdots+(a+12)\phi_1+a+b+1,\\
					&=&\cdots+(a-12)\phi_2-a+b+1,\\
					&=&\cdots-6\phi_2+ax+b+1.
				\end{array}$$
				We conclude that $3$ does not divide $(\mathbb{Z}_K:\mathbb{Z}[\alpha])$ if and only if $\nu_3(a+b+1)=\nu_3(-a+b+1)=\nu_3(ax+b+1)=1$; $(a,b)\in\{(0,2),(0,5),(3,-1),(6,-1)\}\md{9}$.
				
				\item[(c)] If $b\equiv 1\md{3}$, then $F(x)\equiv (x^2-x-1)^3(x^2+x-1)^3\md{3}$. Let $\phi_1(x)=x^2-x-1$ and $\phi_2(x)=x^2+x-1$, then 
				$$\begin{array}{lll}
					F(x)&=&\cdots+(324+420x)\phi_1+(a+144)x+b+89,\\
					&=&\cdots+(324-420x)\phi_2+(a-144)x+b+89.
				\end{array}$$
				We conclude that $3$ does not divide $(\mathbb{Z}_K:\mathbb{Z}[\alpha])$ if and only if $\nu_3((a+144)x+b+89)=\nu_3((a-144)x+b+89)=1$; that is $$(a,b)\in\{(0,4),(0,7),(3,1),(3,4),(3,7),(6,1),(6,4),(6,7)\}\md{9}.$$
			\end{enumerate}
			\item {For $p=11$, assume that $11$ divides $b$}.
			\begin{enumerate}
				\item[(a)] If $11$ divides $a$, then $11$ does not divide $(\mathbb{Z}_K:\mathbb{Z}[\alpha])$ if and only if $\nu_{11}(b)=1$.
				
				\item[(b)] If $11$ does not divide $a$, then $F(x)\equiv x(x+a)^{11}\md{11}$. Let $\phi=x+a$, then $F(x)=\cdots+(a-12a^{11})\phi+b-a^2+a^{12}$. We conclude that $11$ does not divide $(\mathbb{Z}_K:\mathbb{Z}[\alpha])$ if and only if $\nu_{11}(b)=1$ and  $\nu_{11}(b-a^2+a^{12})=1$.
				
			\end{enumerate}
			\item For $p\not\in\{2,3,11\}$ and $p$ divides $a$ and $b$, then $F(x)\equiv x^{12}\md{p}$. Let $\phi=x$, then $F(x)=\phi^{12}+a\phi+b$. Hence, $p$ does not divide $(\mathbb{Z}_K:\mathbb{Z}[\alpha])$ if and only if $\nu_{p}(b)=1$.
			
			\item For $p\not\in\{2,3,11\}$, $p$ does not divide $ab$, and $p^2$ divides $2^{24}\times 3^{12}b^{11}-11^{11}a^{12}$, then by  \cite[Theorem $1$]{LNV}:
			\begin{equation*}
				\nu_p(d_K)=\left\{
				\begin{array}{ll}
					0&\text{ if }\nu_p(2^{24}\times 3^{12}b^{11}-11^{11}a^{12})\text{ is even},\\
					1&\text{ if }\nu_p(2^{24}\times 3^{12}b^{11}-11^{11}a^{12})\text{ is odd}.
				\end{array}
				\right.	
			\end{equation*}
			Thus, by the index formula $(1.1)$, $p$ divides the index $(\mathbb{Z}_K:\mathbb{Z}[\alpha])$.
			
		\end{enumerate}
		\begin{flushright}
			$\square$
		\end{flushright}

		\textit{\textbf{Proof of Theorem \ref{thmp2}}}.\\
		Let $F(x)=x^{12}+ax+b$. Since $\Delta(F)=2^{24}\times 3^{12}b^{11}-11^{11}a^{12}$ and thanks to the index formula $(1.1)$, if $2$ does not divide $a$, then $2$ does not divide $(\mathbb{Z}_{K}:\mathbb{Z}[\alpha])$ and so $2$ does not divide $i(K)$. Now assume that $2$ divides $a$, then we have the following cases:
		\begin{enumerate}
			\item If $\nu_2(b)=0$, then $F(x)\equiv (x-1)^4(x^2+x+1)^4\md{2}$. Let $\phi_1=x-1$ and $\phi_2=x^2+x+1$, then 
			$$
			\begin{array}{lll}
				F(x)&=&\cdots+495\phi_2^4+220\phi_1^3+66\phi_2^2+(a+12)\phi_1+a+b+1,\\
				&=&\cdots+(-25-5x)\phi_2^4+(18+24x)\phi_2^3-18x\phi_1^2+(-4+4x)\phi_1+ax+b+1.
			\end{array}
			$$
			Let $u=\nu_2(a+b+1)$ and $v=\nu_2(a+12)$.
			\begin{enumerate}
				\item[(i)] If $(a,b)\in\{(0,1),(2,3)\}\md{4}$, then by Theorem $\ref{mong}$, $2$ does not divides $i(K)$.
				
				\item[(ii)] If $(a,b)\equiv (0,3)\md{4}$, then for $(a,b)\in\{(0,3),(0,7)\}\md{8}$, by Theorem $\ref{thm1}$, $2$ divides $i(K)$. In particular, if $(a,b)\equiv (0,3)\md{8}$, then $N_{\phi_2}^+(F)=S_2$ has a single side joining $(0,2)$ and $(2,0)$ with $R_{\lambda_2}(F)(y)=(x+1)y^2+xy+1=(y+1)(y+x)\in\mathbb{F}_{\phi_2}[y]$. On the other hand, $N_{\phi_1}^+(F)=S_1$ has a single side joining $(0,2)$ and $(2,0)$ with $R_{\lambda_2}(F)(y)=y^2+y+1$, which is irreducible over $\mathbb{F}_{\phi_1}$. Hence, $2\mathbb{Z}_K=\mathfrak{p}_1^2\mathfrak{p}_{21}^2\mathfrak{p}_{22}^2$ with residue degree $2$ each. By \cite[Theorem 8]{En}, $\nu_2(i(K))=2$.
				
				 If $(a,b)\equiv (4,7)\md{8}$, then $N_{\phi_i}^+(F)=S_i$ has a single side joining $(2,0)$ and $(4,0)$ with $R_{\lambda_1}(F)(y)=y^2+y+1$ and $R_{\lambda_2}(F)(y)=(x+1)y^2+xy+x$, which are irreducible over $\mathbb{F}_{\phi_1}$ and $\mathbb{F}_{\phi_2}$ respectively. Thus $2\mathbb{Z}_{K}=\mathfrak{p}_1^2\mathfrak{p}_2^2$ with $f_1=1$ and $f_2=4$. If $(a,b)\equiv (4,3)\md{8}$, then $N_{\phi_2}^+(F)=S_2$ has a single side joining $(0,2)$ and $(2,0)$ with $R_{\lambda_2}(F)(y)=(x+1)y^2+xy+x+1$, which is irreducible over $\mathbb{F}_{\phi_2}$. Thus $\phi_2$ provides a unique prime ideal of $\mathbb{Z}_{K}$ lying above $2$ with residue degree $4$. On the other hand, we have $u\geq3$ and $v\geq3$. If $u>2v-1$, then $N_{\phi_1}^+(F)=S_{11}+S_{12}+S_{13}$ has three sides of degree $1$ each. Thus $2\mathbb{Z}_{K}=\mathfrak{p}_{11}\mathfrak{p}_{12}\mathfrak{p}_{13}^2\mathfrak{p}_2^4$ with $f_{11}=f_{12}=f_{13}=1$ and $f_2=2$. Hence $2$ divides $i(K)$ and by \cite[Corollary, page 230]{En}, $\nu_2(i(K))=1$. If $u=2v-1$, then $N_{\phi_1}^+(F)=S_{11}+S_{12}$ has two sides joining $(0,u)$, $(2,1)$, and $(4,0)$ with $R_{\lambda_{11}}(F)(y)=y^2+y+1$, which is irreducible over $\mathbb{F}_{\phi_1}$. Thus $2\mathbb{Z}_{K}=\mathfrak{p}_{11}\mathfrak{p}_{12}^2\mathfrak{p}_2^4$ with $f_{11}=2$, $f_{12}=1$, and $f_{2}=2$. Hence $2$ does not divide $i(K)$. If $u<2v-1$, and $u$ is even, then $N_{\phi_1}^+(F)=S_{11}+S_{12}$ has two sides of degree $1$ each. Thus $2\mathbb{Z}_{K}=\mathfrak{p}_{11}^2\mathfrak{p}_{12}^2\mathfrak{p}_2^4$ with $f_{11}=f_{12}=1$ and $f_2=2$. Hence $2$ does not divide $i(K)$. If $u<2v-1$ and $u$ is odd, then $N_{\phi_1}^+(F)=S_{11}+S_{12}$ has two sides joining $(0,u)$, $(2,1)$, and $(4,0)$ with $R_{\lambda_{11}}(F)(y)=(y+1)^2\in\mathbb{F}_{\phi_1}[y]$. Let $s$ be an odd integer such that $F(x)$ is $(x-s)$-regular with respect to $2$, then $F(x)=\cdots+495s^8(x-s)^4+220s^9(x-s)^3+66s^{10}(x-s)^2+(a+12s^{11})(x-s)+b+as+s^{12}$. Since $\phi_2$ provides a unique prime ideal of $\mathbb{Z}_{K}$ lying above $2$ with residue degree $4$, then $2$ divides $i(K)$ if and only if $\nu_2(b+as+s^{12})>2\nu_2(a+12s^{11})-1$. In this case, $2\mathbb{Z}_{K}=\mathfrak{p}_{11}\mathfrak{p}_{12}\mathfrak{p}_{13}^2\mathfrak{p}_2^4$ with $f_{11}=f_{12}=f_{13}=1$ and $f_2=2$. By \cite[Corollary, page 230]{En}, $\nu_2(i(K))=1$.
				
				\item[(iii)] If $(a,b)\equiv(2,1)\md{4}$, then $N_{\phi_2}^+(F)$ has a single side of height $1$. Thus $\phi_2$ provides a unique prime ideal of $\mathbb{Z}_{K}$ lying above $2$ with residue degree $2$. On the other hand, $N_{\phi_1}^+(F)=S_{11}+S_{12}$ has two sides joining $(0,u)$, $(1,1)$, and $(4,0)$ with $u\geq 2$. Thus $2\mathbb{Z}_{K}=\mathfrak{p}_{11}\mathfrak{p}_{12}^3\mathfrak{p}_2^4$ with $f_{11}=f_{12}=1$ and $f_{2}=2$. Hence $2$ does not divide $i(K)$.
			\end{enumerate}
		
			\item If $\nu_2(b)\geq 1$, then $F(x)\equiv x^{12}\md{2}$. Let $\phi=x$, then $F(x)=\phi^{12}+a\phi+b$.\\
			\textbf{If $\nu_2(b)>\cfrac{12}{11}\nu_2(a)$}, then $N_{\phi}(F)$ has two sides of degree $1$ each, because $\nu_2(a)\leq 10$ by assumption. Thus $2\mathbb{Z}_{K}=(1,11)(1,1)$.\\
			\textbf{If $\nu_2(b)<\cfrac{12}{11}\nu_2(a)$}, then we have the following cases:
			\begin{enumerate}
				
				\item[(i)] If gcd$(\nu_2(b),12)=1$, then $N_{\phi}(F)$ has a single side of degree $1$. Thus $2\mathbb{Z}_{K}=\mathfrak{p}^{12}$ with residue degree $1$ and so $2$ does not divide $i(K)$.
				
				\item[(ii)] If $\nu_2(b)\in\{2,10\}$, then $N_{\phi}(F)=S$ has two sides joining $(0,\nu_2(b))$ and $(12,0)$ with $R_{\lambda_1}(F)(y)=(y+1)^2\in \mathbb{F}_{\phi}[y]$. Since $-\lambda=-1/6,-5/6$ is the slope of $S$ respectively, then $6$ divides the ramification index $e$ of any prime ideal of $\mathbb{Z}_{K}$ lying above $2$. Thus the possible cases are: $2\mathbb{Z}_{K}=\mathfrak{p}^6$ with $f=6$, $2\mathbb{Z}_{K}=\mathfrak{p}^{12}$ with $f=1$, and $2\mathbb{Z}_{K}=\mathfrak{p}_1^6\mathfrak{p}_2^6$ with $f_1=f_2=1$. Hence $2$ does not divide $i(K)$.
				
				\item[(iii)] If $\nu_2(b)\in\{3,9\}$, then $N_{\phi}(F)=S$ has a single side joining $(0,\nu_2(b))$ and $(12,0)$ with $R_{\lambda}(F)(y)=y^3+1=(y+1)(y^2+y+1)\in\mathbb{F}_{\phi}[y]$. Thus $2\mathbb{Z}_{K}=\mathfrak{p}_1^4\mathfrak{p}_2^4$ with $f_1=1$ and $f_2=2$. Hence $2$ does not divide $i(K)$. 
				
				\item[(iv)] If $\nu_2(b)=4$; $b\equiv 16\md{32})$ and $a\equiv 0\md{16}$, then $N_{\phi}(F)=S$ has a single side joining $(0,4)$ and $(12,0)$ with $R_{\lambda}(F)(y)=(y+1)^4\in\mathbb{F}_{\phi}[y]$. In this case, we have to use Newton polygon of second order. Let $\omega_2=3[\nu_2,\phi,1/3]$ be the valuation of second order Newton polygon and $g_2=x^3-2$ the key polynomial of $\omega_2$, where $[\nu_2,\phi,1/3]$ is the augmented valuation of $\nu_2$ with respect to $\phi$ and $\lambda=1/3$. Let $F(x)=g_2^4+8g_2^3+24g_2^2+32g_2+ax+b+16$ the $g_2$-expansion of $F(x)$ and $N_2(F)$ the $g_2$-Newton polygon of $F(x)$ with respect to $\omega_2$. We have $\omega_2(x)=1$, $\omega_2(g_2)=3$, and $\omega_2(m)=3\times \nu_2(m)$ for every $m\in\mathbb{Q}_2$.
				
				\begin{enumerate}
					\item[(a)] If $a\equiv 16\md{32}$, then $N_2(F)=T$ has a single side joining $(0,13)$ and $(4,12)$. Thus $2\mathbb{Z}_{K}=\mathfrak{p}^{12}$ with residue degree $1$ and so $2$ does not divide $i(K)$.
					
					\item[(b)] If $b\equiv 16\md{64}$ and $a\equiv0\md{32}$, then $N_2(F)=T$ has a single side joining $(0,15)$ and $(4,12)$. Thus $2\mathbb{Z}_{K}=\mathfrak{p}^{12}$ with residue degree $1$ and so $2$ does not divide $i(K)$.
					
					\item[(c)] If $b\equiv 48\md{64}$ and $a\equiv32\md{64}$, then for $g_2=x^3-2x-2$, we have  $F(x)=g_2^4+(8+8x)g_2^3+(24x^2+48x+56)g_2^2+(96x^2+176x+160)g_2+128x^2+(a+224)x+b+144$. If $(a,b)\in\{(48,32),(48,96)\}\md{128}$, then $\omega_2(128x^2+(a+224)x+b+144)=18$, if $(a,b)\equiv (112,32)\md{128}$, then $\omega_2(128x^2+(a+224)x+b+144)\geq21$, and finally if $(a,b)\equiv (112,96)\md{128}$, then  $\omega_2(128x^2+(a+224)x+b+144)=19$. In every case, $N_2(F)=T_1+T_2$ has two sides joining $(0,w)$, $(1,16)$, and $(4,12)$ with $w\geq 18$. Thus $2\mathbb{Z}_{K}=\mathfrak{p}_1^3\mathfrak{p}_2^9$ with residue degree $1$ each and so $2$ does not divide $i(K)$.
					
					\item[(d)] If $b\equiv 48\md{128}$ and $a\equiv0\md{64}$, then $N_2(F)=T$ has a single side joining $(0,18)$ and $(4,12)$ with $R_2(y)=y^2+y+1$, which is irreducible over $\mathbb{F}_2$. Thus $2\mathbb{Z}_{K}=\mathfrak{p}^6$ with residue degree $2$. Hence $2$ does not divide $i(K)$.
					
					\item[(e)] If $b\equiv 112\md{256}$ and $a\equiv64\md{128}$, then for $g_2=x^3-2x^2-2$, we have $F(x)=g_2^4+(8x^2+24x+88)g_2^3+(288x^2+336x+952)g_2^2+(2080x^2+1312x+3296)g_2+4096x^2+(a+1472)x+b+3472$. If $(a,b)\equiv(64,112)\md{256}$, then $\omega_2(4096x^2+(a+1472)x+b+3472)\geq24$, and if $(a,b)\equiv(192,112)\md{256}$, then $\omega_2(4096x^2+(a+1472)x+b+3472=22$. Thus $N_2(F)=T_1+T_2+T_3$ has three sides joining $(0,w)$, $(1,18)$, $(2,15)$, and $(4,12)$ with $w\geq22$. Hence, $2\mathbb{Z}_{K}=\mathfrak{p}_1^3\mathfrak{p}_2^3\mathfrak{p}_3^6$ with residue degree $1$ each. Hence $2$ divides $i(K)$.
					
					\item[(f)] If $b\equiv 112\md{256}$ and $a\equiv0\md{128}$, then $N_2(F)=T_1+T_2$ has two sides joining $(0,21)$, $(2,15)$, and $(4,12)$ with $R_{21}(y)=y^2+y+1$, which is irreducible over $\mathbb{F}_2$. Thus $2\mathbb{Z}_{K}=\mathfrak{p}_1^3\mathfrak{p}_2^6$ with $f_1=2$ and $f_2=1$. Hence $2$ does not divide $i(K)$.
					
					\item[(g)] If $b\equiv 240\md{256}$ and $a\equiv0\md{128}$, then $N_2(F)=T_1+T_2+T_3$ has three sides joining $(0,w)$, $(1,18)$, $(2,15)$, and $(4,12)$ with $w\geq22$. Thus $2\mathbb{Z}_{K}=\mathfrak{p}_1^3\mathfrak{p}_2^3\mathfrak{p}_3^6$ with residue degree $1$ each. Hence $2$ divides $i(K)$.
				\end{enumerate}

				\item[(v)] If $\nu_2(b)=6$; $b\equiv 64\md{128})$ and $a\equiv 0\md{64}$, then $N_{\phi}(F)=S_1$ has a single side joining $(0,6)$ and $(12,0)$ with $R_{\lambda_1}(F)(y)=y^6+1=(y+1)^2(y^2+y+1)^2\in\mathbb{F}_{\phi}[y]$. In this case, we have to use Newton polygon of second order. Let $\omega_2=2[\nu_2,\phi,1/2]$ be the valuation of second order Newton polygon. Let $g_{2,1}=x^2-2$ and $g_{2,2}=x^4-2x^2-4$ be the key polynomials of $(\omega_2,y+1)$ and $(\omega_2,y^2+y+1)$ respectively, then 
				$$
				\begin{array}{lll}
					F(x)&=&\cdots+240g_{2,1}^2+192g_{2,1}+ax+b+64,\\
					&=&g_{2,2}^3+(24+6x^2)g_{2,2}^2+(160+80x^2)g_{2,2}+256x^2+ax+b+320.
				\end{array}
				$$
				We have $\omega_2(x)=1$, $\omega_2(g_{2,1})=2$, $\omega_2(g_{2,2})=4$, and $\omega_2(m)=2\times \nu_2(m)$ for every $m\in\mathbb{Q}_2$.
				\begin{enumerate}
					\item[(a)] If $a\equiv 64\md{128}$, then $N_{2,i}^+(F)=T_{2i}$ has a single side joining $(0,13)$ and $(2,12)$ for every $i=1,2$. Thus $2\mathbb{Z}_{K}=\mathfrak{p}_1^4\mathfrak{p}_2^4$ with $f_1=1$ and $f_2=2$. Hence $2$ does not divide $i(K)$.
					
					\item[(b)] If $b\equiv 64\md{256}$ and $a\equiv 0\md{128}$, then for $g_{2,1}=x^2-2x-2$ and $g_{2,2}=x^4-2x^3-2x^2-4x-12$, let  
					$$
					\begin{array}{lll}
						F(x)&=&\cdots+(31536+23040x)g_{2,1}^2+(54656+54976x)g_{2,1}+(a+49920)x+b+36544,\\
						&=&g_{2,2}^3+(6x^3+18x^2+68x+272)g_{2,2}^2+(992x^3+1728x^2+4320x+11456)g_{2,2}\\
						&&+27520x^3+35840x^2+(a+63488)x+b+100032.
					\end{array}
					$$
					We have $\omega_2((54656+54976x)g_{2,1})=\omega_2((992x^3+1728x^2+4320x+11456)g_{2,2})=15$.
					
					\item[(b$_1$)] If $a\equiv 128\md{256}$, then $N_{2,i}^+(F)=T_{i}$ has a single side joining $(0,15)$ and $(2,12)$ for every $i=1,2$. Thus $2\mathbb{Z}_{K}=\mathfrak{p}_1^4\mathfrak{p}_2^4$ with $f_1=1$ and $f_2=2$. Hence $2$ does not divide $i(K)$.
					
					\item[(b$_2$)] If $(a,b)\equiv (0,320)\md{512}$, then $N_{2,i}^+(F)=T_{i}$ has a single side  joining $(0,17)$ and $(2,12)$ for $i=1,2$. Hence, $2\mathbb{Z}_{K}=\mathfrak{p}_1^4\mathfrak{p}_2^4$ with $f_1=1$ and $f_2=2$. Hence $2$ does not divide $i(K)$.
					
					\item[(b$_3$)] If $(a,b)\equiv (256,320)\md{512}$, then $\omega_2(27520x^3+35840x^2+(a+63488)x+b+100032)=17$. Thus $N_{2,2}^+(F)=T_{2}$ has a single side  joining $(0,17)$ and $(2,12)$. If $(a,b)\in\{(256,320),(768,230)\}\md{1024}$, then 	$N_{2,1}^+(F)=T_{11}+T_{12}$ has two sides joining $(0,w)$, $(1,15)$, and $(2,12)$ with $w\geq 19$, and so $2\mathbb{Z}_{K}=\mathfrak{p}_{11}^2\mathfrak{p}_{12}^2\mathfrak{p}_2^4$ with $f_{11}=f_{12}=1$ and $f_2=2$. Hence $2$ does not divide $i(K)$. If $(a,b)\in\{(256,832),(768,832)\}\md{1024}$, then $N_{2,1}^+(F)=T_{1}$ has a single side  joining $(0,18)$ and $(2,12)$ with $R_{21}(y)=y^2+y+1$, which is irreducible over $\mathbb{F}_{21}$. Thus $2\mathbb{Z}_{K}=\mathfrak{p}_1^2\mathfrak{p}_2^4$ with residue degree $2$ each. Hence $2$ divides $i(K)$.
					
					\item[(b$_4$)] If $(a,b)\in\{(0,64),(256,64)\}\md{512}$, then for $g_{2,1}=x^2-2x-6$ and $g_{2,2}=x^4-2x^3-6x^2-4x-12$, let 
					
					$$\begin{array}{lll}
						F(x)&=&\cdots+(304560+123264x)g_{2,1}^2+(1141504+573632x)g_{2,1}+(a+1041920)x+b+1715136,\\
						&=&g_{2,2}^3+(6x^3+30x^2+116x+560)g_{2,2}^2+(2496x^3+6368x^2+12768x+43264)g_{2,2}\\
						&&+145152x^3+285184x^2+(a+257536)x+b+440256.
					\end{array}$$ 
					We have $\omega_2((1141504+573632x)g_{2,1})=\omega_2((+2496x^3+6368x^2+12768x+43264)g_{2,2})=15$. If $a\equiv 256\md{512}$, then $N_{2,i}^+(F)=T_{i}$ has a single side  joining $(0,17)$ and $(2,12)$ for $i=1,2$. Thus $2\mathbb{Z}_{K}=\mathfrak{p}_1^4\mathfrak{p}_2^4$ with $f_1=1$ and $f_2=2$. Hence $2$ does not divide $i(K)$. If $(a,b)\in\{(0,576),(512,576)\}\md{1024}$, then $N_{2,i}^+(F)=T_{i}$ has a single side  joining $(0,18)$ and $(2,12)$ with $R_{21}(y)=y^2+y+1$ and $R_{22}(y)=y^2+xy+1$, which are irreducible over $\mathbb{F}_{21}$ and  $\mathbb{F}_{22}$ respectively. Thus $2\mathbb{Z}_{K}=\mathfrak{p}_1^2\mathfrak{p}_2^2$ with $f_1=2$ and $f_2=4$. Hence $2$ does not divide $i(K)$. If $(a,b)\in\{(0,64),(512,64)\}\md{1024}$, then $N_{2,i}^+(F)=T_{i1}+T_{i2}$ has two sides  joining $(0,w_i)$, $(1,15)$, and $(2,12)$ with $w_1\geq19$ and $w_2=19$. Thus $2\mathbb{Z}_{K}=\mathfrak{p}_{11}^2\mathfrak{p}_{12}^2\mathfrak{p}_{21}^2\mathfrak{p}_{22}^2$ with $f_{11}=f_{12}=1$ and $f_{21}=f_{22}=2$. Hence $2$ divides $i(K)$ and by \cite[Theorem 7]{En}, $\nu_2(i(K))=3$.
					
					\item[(c)] If $b\equiv 192\md{256}$ and $a\equiv 0\md{128}$, then $\omega_2(ax+b+64)\geq 15$ and $\omega_2(256x^2+ax+b+320)\geq 15$.
					
					\item[(c$_1$)] If $a\equiv 128\md{256}$, then $N_{2,i}^+(F)=T_{i}$ has a single side  joining $(0,15)$ and $(2,12)$ for $i=1,2$. Thus $2\mathbb{Z}_{K}=\mathfrak{p}_1^4\mathfrak{p}_2^4$ with $f_1=1$ and $f_2=2$. Hence $2$ does not divide $i(K)$.
					
					\item[(c$_2$)] If $(a,b)\in\{(0,448),(256,448)\}\md{512}$, $N_{2,1}^+(F)=T_{11}+T_{12}$ has two sides  joining $(0,w)$, $(1,14)$, and $(2,12)$ with $w\geq17$. On the other hand we have $N_{2,2}^+(F)=T_{2}$ has a single side joining $(0,16)$ and $(2,12)$ with $R_{22}(y)=y^2+xy+1$, which is irreducible over $\mathbb{F}_{22}$. Thus $2\mathbb{Z}_{K}=\mathfrak{p}_{11}^2\mathfrak{p}_{12}^2\mathfrak{p}_{2}^2$ with $f_{11}=f_{12}=1$ and $f_2=4$. Hence $2$ does not divide $i(K)$.
					
					\item[(c$_3$)] If $(a,b)\in\{(0,192),(256,192)\}\md{512}$, then $N_{2,1}^+(F)=T_{1}$ has a single side  joining $(0,16)$ and $(2,12)$ with $R_{21}(y)=y^2+y+1$, which is irreducible over $\mathbb{F}_{21}$. On the other hand, $N_{2,2}^+(F)=T_{21}+T_{22}$ has two sides joining $(0,u)$, $(1,14)$, and $(2,12)$ with $u\geq17$. Thus $2\mathbb{Z}_{K}=\mathfrak{p}_1^2\mathfrak{p}_{21}^2\mathfrak{p}_{22}^2$ with residue degree $2$ each. Hence $2$ divides $i(K)$ and by \cite[Theorem 8]{En}, $\nu_2(i(K))=2$. 
				\end{enumerate}
				\item[(vi)] If $\nu_2(b)=8$; $b\equiv 256\md{512})$ and $a\equiv 0\md{256}$, then $N_{\phi}(F)=S$ has a single side joining $(0,8)$ and $(12,0)$ with $R_{\lambda_1}(F)(y)=(y+1)^4\in\mathbb{F}_{\phi}[y]$. In this case, we have to use Newton polygon of second order. Let $\omega_2=3[\nu_2,\phi,2/3]$ be the valuation of second order Newton polygon and $g_2=x^3-16$ the key polynomial of $\omega_2$. Let $F(x)=g_2^4+16g_2^3+96g_2^2+256g_2+ax+b+256$, then we have  $\omega_2(x)=2$, $\omega_2(g_2)=6$, and $\omega_2(m)=3\times\nu_2(m)$ for every $m\in\mathbb{Q}_2$.
				\begin{enumerate}
					\item[(a)] If $a\equiv 256\md{512}$, then $N_2(F)=T$ has a single side  joining $(0,26)$ and $(4,24)$ with $R_2(y)=(y+1)^2\in\mathbb{F}_2[y]$. Since $-\lambda'=-1/2$ is the slope of $T$ and  $-\lambda=-2/3$ is the slope of $S$, then $6$ divides the ramification index $e$ of any prime ideal of $\mathbb{Z}_{K}$ lying above $2$. Thus the possible cases are: $2\mathbb{Z}_{K}=\mathfrak{p}^{12}$ with $f_=1$, $2\mathbb{Z}_{K}=\mathfrak{p}_1^6\mathfrak{p}_2^6$ with $f_1=f_2=1$, and $2\mathbb{Z}_{K}=\mathfrak{p}^6$ with $f=2$. Hence $2$ does not divide $i(K)$.
					
					\item[(b)] If $b\equiv256\md{1024}$ and $a\equiv 0\md{512}$, then $N_2(F)=T$ has a single side joining $(0,27)$ and $(4,24)$. Thus $2\mathbb{Z}_{K}=\mathfrak{p}^{12}$ with residue degree $1$. Hence $2$ does not divide $i(K)$. 
					
					\item[(c)] If $b\equiv768\md{1024}$ and $a\equiv 512\md{1024}$, then $N_2(F)=T$ has a single side joining $(0,29)$ and $(4,24)$. Thus $2\mathbb{Z}_{K}=\mathfrak{p}^{12}$ with residue degree $1$. Hence $2$ does not divide $i(K)$. 
					
					\item[(d)] If $b\equiv768\md{2^{11}}$ and $a\equiv 0\md{1024}$, then $N_2(F)=T$ has a single side joining $(0,30)$ and $(4,24)$ with $R_2(y)=y^2+y+1$, which is irreducible over $\mathbb{F}_2$. Thus $2\mathbb{Z}_{K}=\mathfrak{p}^{6}$ with residue degree $2$. Hence $2$ does not divide $i(K)$. 
					
					\item[(e)] If $b\equiv768+2^{10}\md{2^{11}}$ and $a\equiv 2^{10}\md{2^{11}}$, then $N_2(F)=T_{1}+T_{2}$ has two sides joining $(0,32)$, $(2,27)$, and $(4,24)$. Thus $2\mathbb{Z}_{K}=\mathfrak{p}_1^6\mathfrak{p}_2^6$ with residue degree $1$ each. Hence $2$ does not divide $i(K)$. 
					
					\item[(f)] $b\equiv768+2^{10}\md{2^{12}}$ and $a\equiv 0\md{2^{11}}$, then $N_2(F)=T_{1}+T_{2}$ has two sides joining $(0,33)$, $(2,27)$, and $(4,24)$ with $R_{21}(y)=y^2+y+1$, which is irreducible over $\mathbb{F}_2$. Thus $2\mathbb{Z}_{K}=\mathfrak{p}_1^3\mathfrak{p}_2^6$ with $f_1=2$ and $f_2=1$. Hence $2$ does not divide $i(K)$. 
					
					\item[(g)] If $b\equiv 768+2^{10}+2^{11}\md{2^{12}}$ and $a\equiv 0\md{2^{11}}$, then $N_2(F)=T_{1}+T_{2}+T_{3}$ has three sides joining $(0,w)$, $(1,30)$, $(2,27)$, and $(4,24)$ with $w\geq34$. Hence, $2\mathbb{Z}_{K}=\mathfrak{p}_1^3\mathfrak{p}_2^3\mathfrak{p}_3^6$ with residue degree $1$ each. Hence $2$ divides $i(K)$.
					
				\end{enumerate}
				
			\end{enumerate}
		\end{enumerate}
		\begin{flushright}
			$\square$
		\end{flushright}
		
		\textit{\textbf{Proof of Theorem \ref{thmp3}}}.\\
		Since $\Delta(F)=2^{24}\times 3^{12}b^{11}-11^{11}a^{12}$ and thanks to the index formula $(1.1)$, if $3$ does not divide $a$, then $3$ does not divide $(\mathbb{Z}_{K}:\mathbb{Z}[\alpha])$, and so $3$ does not divide $i(K)$. Now assume that $3$ divides $a$, then we have the following cases:
		\begin{enumerate}
			\item If $\nu_3(b)=0$ and $b\equiv -1\md{3}$, then $F(x)\equiv(x-1)^3(x+1)^3(x^2+1)^3\md{3}$. Let $\phi_1=x-1$, $\phi_2=x+1$, and $\phi_3=x^2+1$, then 
			
			$$
			\begin{array}{lll}
				F(x)&=&\cdots+220\phi_1^3+66\phi_1^2+(a+12)\phi_1+b+a+1,\\
				&=&\cdots-220\phi_2^3+66\phi_2^2+(a-12)\phi_2+b-a+1,\\
				&=&\cdots-\phi_3^3+15\phi_3^2-6\phi_3+ax+b+1.
			\end{array}
			$$
			\begin{enumerate}
				\item[(i)] If $(a,b)\in\{(0,2),(0,5),(3,-1),(6,-1)\}\md{9}$, then by Theorem $\ref{mong}$, $3$ does not divide $i(K)$.
				
				\item[(ii)] If $(a,b)\equiv (0,-1)\md{9}$, then by Theorem $\ref{thm1}$, $3$ divides $i(K)$.
				
				\item[(iii)] If $(a,b)\in\{(3,5),(6,2)\}\md{9}$, then $N_{\phi_i}^+(F)$ has a single side of height $1$ for every $i=2,3$. Thus there are two prime ideals $\mathfrak{p}_2^3$ and $\mathfrak{p}_3^2$ of $\mathbb{Z}_K$ lying above $3$ with $f_2=1$ and $f_3=2$, provided by $\phi_2$ and $\phi_3$ respectively. Let $u=\nu_3(a+b+1)$ and $v=\nu_3(a+12)$.
				
				\begin{enumerate}
					\item[(a)] If $(a,b)\equiv (3,5)\md{9}$, then $N_{\phi_1}^+(F)=S_{11}+S_{12}$ has two sides joining $(0,u)$, $(1,1)$, and $(3,0)$ with $u\geq2$. Thus $3\mathbb{Z}_{K}=\mathfrak{p}_{11}\mathfrak{p}_{12}^2\mathfrak{p}_2^3\mathfrak{p}_3^3$ with $f_{11}=f_{12}=f_2=1$ and $f_3=2$. Hence $2$ does not divide $i(K)$.
					
					\item[(b)] If $(a,b)\in\{(6,2),(6,11),(15,2),(15,20),(24,11),(24,20)\}\md{27}$, then $N_{\phi_1}^+(F)=S_1$ has a single side joining $(0,2)$ and $(3,0)$. Thus $3\mathbb{Z}_{K}=\mathfrak{p}_{1}^3\mathfrak{p}_2^3\mathfrak{p}_3^3$ with $f_1=f_2=1$ and $f_3=2$. Hence $2$ does not divide $i(K)$.
					
					\item[(c)] If $(a,b)\in\{(15,11),(15,38),(42,11),(42,65),(69,38),(69,65)\}\md{81}$, then $N_{\phi_1}^+(F)=S_1$ has a single side joining $(0,3)$ and $(3,0)$ with $R_{\lambda_1}(F)(y)=y^3+y^2-1$ or $R_{\lambda_1}(F)(y)=y^3+y^2+1=(y+1)(y^2-y-1)$. Thus $3\mathbb{Z}_{K}=\mathfrak{p}_1\mathfrak{p}_2^3\mathfrak{p}_3^3$ with $f_1=3$, $f_2=1$, and $f_3=2$ or $3\mathbb{Z}_{K}=\mathfrak{p}_{11}\mathfrak{p}_{12}\mathfrak{p}_2^3\mathfrak{p}_3^3$ with $f_{11}=f_2=1$ and $f_{12}=f_3=2$ respectively. Hence, $3$ does not divide $i(K)$.
					
					\item[(d)] If $(a,b)\in\{(15,65),(42,38),(69,11)\}\md{81}$, then $u\geq4$ and $v\geq3$.
					
					\item[(d$_1$)] If $u>2v-1$, then $N_{\phi_1}^+(F)$ has three sides with degree $1$ each. Thus $3\mathbb{Z}_{K}=\mathfrak{p}_{11}\mathfrak{p}_{12}\mathfrak{p}_{13}\mathfrak{p}_2^3\mathfrak{p}_3^3$ with $f_{11}=f_{12}=f_{13}=f_2=1$ and $f_3=2$. Hence, $3$ divides $i(K)$ and by \cite[Corollary, page 230]{En}, $\nu_3(i(K))=1$. 
					
					\item[(d$_2$)] If $u<2v-1$, then $N_{\phi_2}^+(F)=S_{11}+S_{12}$ has two sides joining $(0,u)$, $(2,1)$, and $(3,0)$ with $d(S_{12})=1$. If $u$ is even, then $d(S_{11})=1$. Thus $3\mathbb{Z}_{K}=\mathfrak{p}{11}^2\mathfrak{p}_{12}\mathfrak{p}_2^3\mathfrak{p}_3^3$ with $f_{11}=f_{12}=f_2=1$ and $f_3=2$. Hence $3$ does not divide $i(K)$. If $u$ is odd, then $d(S_{11})=2$ with $R_{\lambda_{11}}(F)(y)=y^2+(a+b+1)_3\in\mathbb{F}_{\phi_1}[y]$. Thus, if $(a+b+1)_3\equiv 1\md{3}$, then $R_{\lambda_{11}}(F)(y)$ is irreducible and so $3\mathbb{Z}_{K}=\mathfrak{p}{11}\mathfrak{p}_{12}\mathfrak{p}_2^3\mathfrak{p}_3^3$ with $f_{11}=f_3=2$ and $f_{12}=f_2=1$. Hence $3$ does not divide $i(K)$. If $(a+b+1)_3\equiv -1\md{3}$, then $R_{\lambda_{11}}(F)(y)=(y-1)(y+1)$. Thus $3\mathbb{Z}_{K}=\mathfrak{p}{111}\mathfrak{p}{112}\mathfrak{p}_{12}\mathfrak{p}_2^3\mathfrak{p}_3^3$ with $f_{111}=f_{112}=f_{12}=f_2=1$ and $f_3=2$. Hence $3$ divides $i(K)$ and by \cite[Corollary, page 230]{En}, $\nu_3(i(K))=1$.

					\item[(d$_3$)] If $u=2v-1$, then $N_{\phi_1}^+(F)=S_{11}+S_{12}$ has two sides joining $(0,u)$, $(2,1)$, and $(3,0)$ with $d(S_{12})=1$ and $R_{\lambda_{11}}(F)(y)=y^2+(a+12)_3y+(a+b+1)_3\in\mathbb{F}_{\phi_1}[y]$. If $(a+b+1)_3\equiv -1\md{3}$, then $R_{\lambda_{11}}(F)(y)$ is irreducible. Thus $3\mathbb{Z}_{K}=\mathfrak{p}_{11}\mathfrak{p}_{12}\mathfrak{p}{2}^3\mathfrak{p}_3^3$ with $f_{11}=f_3=2$ and $f_{12}=f_2=1$. Hence $3$ does not divide $i(K)$. If $(a+b+1)_3\equiv 1\md{3}$, then $R_{\lambda_{21}}(F)(y)=(y\pm1)^2\in\mathbb{F}_{\phi_2}[y]$. Let $s$ be an integer such that $F(x)$ is $(x-s)$-regular with respect to $3$ and $s\equiv 1\md{3}$, then $F(x)=220s^9(x-s)^3+66s^{10}(x-s)^2+(a+12s^{11})(x-s)+b+as+s^{12}$. Since $\nu_3(66s^{10})=1$, then $3$ divides $i(K)$ if and only if  $\nu_3(b+as+s^{12})>2\nu_3(a+12s^{11})-1$ or $\nu_3(b+as+s^{12})<2\nu_3(a+12s^{11})-1$, $\nu_3(b+as+s^{12})$ is odd, and $(b+as+s^{12})_3\equiv -1\md{3}$. In this case, $3\mathbb{Z}_{K}=3\mathbb{Z}_{K}=\mathfrak{p}{111}\mathfrak{p}{112}\mathfrak{p}_{12}\mathfrak{p}_2^3\mathfrak{p}_3^3$ with $f_{111}=f_{112}=f_{12}=f_2=1$ and $f_3=2$. Hence, by \cite[Corollary, page 230]{En}, $\nu_3(i(K))=1$.
					
					\item[(e)] If $(a,b)\in\{(6,20),(33,74),(60,47)\}\md{81}$, then $N_{\phi_1}^+(F)=S_1$ has a single side joining $(0,3)$ and $(3,0)$ with $R_{\lambda_1}(F)(y)=y^3+y^2-y+1$, which is irreducible over $\mathbb{F}_{\phi_1}$. Thus $3\mathbb{Z}_{K}=\mathfrak{p}_1\mathfrak{p}_2^3\mathfrak{p}_3^3$ with $f_1=3$, $f_2=1$, and $f_3=2$.
					
					\item[(f)] If $(a,b)\in\{(6,74),(33,47),(60,20)\}\md{81}$, then $N_{\phi_1}^+(F)=S_{11}+S_{12}$ has two sides joining $(0,u)$, $(1,2)$, and $(3,0)$ with $u\geq4$. Thus $d(S_{11})=1$ and $R_{\lambda_{12}} (F)(y)=y^2+y-1$, which is irreducible over $\mathbb{F}_{\phi_1}$. Hence, $3\mathbb{Z}_{K}=\mathfrak{p}_{11}\mathfrak{p}_{12}\mathfrak{p}_2^3\mathfrak{p}_3^3$ with $f_{11}=f_2=1$ and $f_{12}=f_3=2$. Hence $2$ does not divide $i(K)$.
					
					\item[(g)] If $(a,b)\in\{(6,47),(33,20),(60,74)\}\md{81}$, then $N_{\phi_1}^+(F)=S_1$ has a single side  joining $(0,3)$ and $(3,0)$ with $R_{\lambda_1}(F)(y)=y^3+y^2-y-1=(y-1)(y+1)^2\in\mathbb{F}_{\phi_1}[y]$. Let $s$ be an integer such that $F(x)$ is $(x-s)$-regular with respect to $3$ and $s\equiv 1\md{3}$, then $F(x)=220s^9(x-s)^3+66s^{10}(x-s)^2+(a+12s^{11})(x-s)+b+as+s^{12}$. Since $\nu_3(66s^{10})=1$, then $3$ divides $i(K)$ if and only if  $\nu_3(b+as+s^{12})>2\nu_3(a+12s^{11})-1$ or $\nu_3(b+as+s^{12})<2\nu_3(a+12s^{11})-1$, $\nu_3(b+as+s^{12})$ is odd, and $(b+as+s^{12})_3\equiv -1\md{3}$. In this case, $3\mathbb{Z}_K=\mathfrak{p}{111}\mathfrak{p}{112}\mathfrak{p}_{12}\mathfrak{p}_2^3\mathfrak{p}_3^3$ with $f_{111}=f_{112}=f_{12}=f_2=1$ and $f_3=2$. Hence, by \cite[Corollary, page 230]{En}, $\nu_3(i(K))=1$.
					
					\item[(h)] If $(a,b)\in\{(24,2),(24,29),(51,2),(51,56),(78,29),(78,56)\}\md{81}$, then $N_{\phi_1}^+(F)=S_1$ has a single side joining $(0,3)$ and $(3,0)$ with $R_{\lambda_1}(F)(y)=y^3+y^2+y-1$ or $R_{\lambda_1}(F)(y)=(y+1)(y^2+1)$. Thus $3\mathbb{Z}_{K}=\mathfrak{p}_1\mathfrak{p}_2^3\mathfrak{p}_3^3(3,1)(2,3)(1,3)$ with $f_1=3$, $f_2=1$, and $f_3=2$ or  $3\mathbb{Z}_{K}=\mathfrak{p}_{11}\mathfrak{p}_{12}\mathfrak{p}_2^3\mathfrak{p}_3^2$ with $f_{11}=f_2=1$ and $f_{12}=f_3=2$ respectively. Hence, $3$ does not divide $i(K)$.
					
					\item[(i)] If $(a,b)\in\{(24,56),(51,29),(78,2)\}\md{81}$, then $N_{\phi_1}^+(F)=S_{11}+S_{12}$ has two sides joining $(0,u)$, $(1,2)$, and $(3,0)$ with $u\geq 4$,  and $R_{\lambda_{12}}(F)(y)=y^2+y+1=(y-1)^2\in\mathbb{F}_{\phi_1}[y]$. Let $s$ be an integer such that $F(x)$ is $(x-s)$-regular with respect to $3$ and $s\equiv 1\md{3}$, then $F(x)=220s^9(x-s)^3+66s^{10}(x-s)^2+(a+12s^{11})(x-s)+b+as+s^{12}$. Since $\nu_3(66s^{10})=1$, then $3$ divides $i(K)$ if and only if  $\nu_3(b+as+s^{12})>2\nu_3(a+12s^{11})-1$ or $\nu_3(b+as+s^{12})<2\nu_3(a+12s^{11})-1$, $\nu_3(b+as+s^{12})$ is odd, and $(b+as+s^{12})_3\equiv -1\md{3}$. In this case, $3\mathbb{Z}_{K}=\mathfrak{p}{111}\mathfrak{p}{112}\mathfrak{p}_{12}\mathfrak{p}_2^3\mathfrak{p}_3^3$ with $f_{111}=f_{112}=f_{12}=f_2=1$ and $f_3=2$. Hence, by \cite[Corollary, page 230]{En}, $\nu_3(i(K))=1$.
				\end{enumerate}

				\item[(iv)] If $(a,b)\in\{(3,2),(6,5)\}\md{9}$, then $N_{\phi_i}^+(F)$ has a single side of height $1$ for every $i=1,3$. Thus there are two prime ideals $\mathfrak{p}_1^3$ and $\mathfrak{p}_3^3$ of $\mathbb{Z}_K$ lying above $3$ with $f_1=1$ and $f_3=2$, provided by $\phi_1$ and $\phi_3$ respectively. Let $\mu=\nu_3(b-a+1)$ and $\tau=\nu_3(a-12)$.
				\begin{enumerate}
					\item[(a)] If $(a,b)\equiv (6,5)\md{9}$, then $N_{\phi_2}^+(F)=S_{21}+S_{22}$ has two sides joining $(0,\mu)$, $(1,1)$, and $(3,0)$ with $\mu\geq2$. Thus $3\mathbb{Z}_{K}=\mathfrak{p}_1^3\mathfrak{p}_{21}\mathfrak{p}_{22}^2\mathfrak{p}_3^3$ with $f_1=f_{21}=f_{22}=1$ and $f_3=2$. Hence $3$ does not divide $i(K)$.
					
					\item[(b)] If $(a,b)\in\{(3,11),(3,20),(12,2),(12,20),(21,2),(21,20)\}\md{27}$, then $N_{\phi_2}^+(F)=S_2$ has a single side  joining $(2,0)$ and $(3,0)$. Thus $3\mathbb{Z}_{K}=\mathfrak{p}_1^3\mathfrak{p}_2^3\mathfrak{p}_3^3$ with $f_1=f_2=1$ and $f_3=2$. Hence $3$ does not divide $i(K)$.
					
					\item[(c)] If $(a,b)\in\{(12,38),(12,65),(39,11),(39,65),(66,11),(66,38)\}\md{81}$, then $N_{\phi_2}^+(F)=S_2$ has a single side joining $(0,3)$ and $(3,0)$ with $R_{\lambda_2}(F)(y)=-y^3+y^2-1$ or $R_{\lambda_2}(F)(y)=-y^3+y^2+1=-(y+1)(y^2+y+2)$. Thus $3\mathbb{Z}_{K}=\mathfrak{p}_1^3\mathfrak{p}_2\mathfrak{p}_3^2$ with $f_1=1$, $f_2=3$, and $f_3=2$ or $3\mathbb{Z}_{K}=\mathfrak{p}_1^3\mathfrak{p}_{21}\mathfrak{p}_{22}\mathfrak{p}_3^2$ with $f_1=f_{21}=1$ and $f_{22}=f_3=2$ respectively. Hence $3$ does not divide $i(K)$.
					
					\item[(d)] If $(a,b)\in\{(12,11),(39,38),(66,65)\}\md{81}$, then  $\mu\geq 4$ and $\tau\geq 3$.
					
					\item[(d$_1$)] If $\mu>2\tau-1$, then $N_{\phi_2}^+(F)$ has three sides with degree $1$ each. Thus $3\mathbb{Z}_{K}=\mathfrak{p}_1^3\mathfrak{p}_{21}\mathfrak{p}_{22}\mathfrak{p}_{23}\mathfrak{p}_3^3$ with $f_1=f_{21}=f_{22}=f_{23}=1$ and $f_3=2$. Hence $3$ divides $i(K)$ and by \cite[Corollary, page 230]{En}, $\nu_3(i(K))=1$. 
					
					\item[(d$_2$)] If $\mu<2\tau-1$, then $N_{\phi_2}^+(F)=S_{21}+S_{22}$ has two sides joining $(0,\mu)$, $(2,1)$, and $(3,0)$ with $d(S_{22})=1$.\\
					If $\mu$ is even, then $d(S_{21})=1$. Thus $3\mathbb{Z}_{K}=\mathfrak{p}_1^3\mathfrak{p}_{21}^2\mathfrak{p}_{22}\mathfrak{p}_{3}^3$ with $f_1=f_{21}=f_{22}=1$ and $f_{3}=2$. Hence $3$ does not divide $i(K)$.\\
					If $\mu$ is odd, then $d(S_{21})=2$ with $R_{\lambda_{21}}(F)(y)=y^2+(b-a+1)_3\in\mathbb{F}_{\phi_2}[y]$. Thus, if $(b-a+1)_3\equiv 1\md{3}$, then $R_{\lambda_{21}}(F)(y)$ is irreducible and so $3\mathbb{Z}_{K}=\mathfrak{p}_1^3\mathfrak{p}_{21}\mathfrak{p}_{22}\mathfrak{p}_{3}^3$ with $f_1=f_{22}=1$ and $f_{21}=f_{3}=2$. Hence $3$ does not divide $i(K)$. If $(b-a+1)_3\equiv -1\md{3}$, then $R_{\lambda_{21}}(F)(y)=(y-1)(y+1)$. Hence, $3\mathbb{Z}_{K}=\mathfrak{p}_1^3\mathfrak{p}_{211}\mathfrak{p}_{212}\mathfrak{p}_{22}\mathfrak{p}_3^3$ with $f_1=f_{211}=f_{212}=f_{22}=1$ and $f_3=2$. Hence $3$ divides $i(K)$ and by \cite[Corollary, page 230]{En}, $\nu_3(i(K))=1$. 
					
					\item[(d$_3$)] If $\mu=2\tau-1$, then $N_{\phi_2}^+(F)=S_{21}+S_{22}$ has two sides joining $(0,\mu)$, $(2,1)$, and $(3,0)$ with $d(S_{22})=1$ and $R_{\lambda_{21}}(F)(y)=y^2+(a-12)_3y+(a-b+1)_3\in\mathbb{F}_{\phi_2}[y]$. If $(b-a+1)_3\equiv -1\md{3}$, then $R_{\lambda_{21}}(F)(y)$ is irreducible. Thus $3\mathbb{Z}_{K}=\mathfrak{p}_1^3\mathfrak{p}_{21}\mathfrak{p}_{22}\mathfrak{p}_3^3$ with $f_1=f_{22}=1$ and $f_{22}=f_3=2$. Hence $3$ does not divide $i(K)$. If $(b-a+1)_3\equiv 1\md{3}$, then $R_{\lambda_{21}}(F)(y)=(y\pm1)^2\in\mathbb{F}_{\phi_2}[y]$. Let $s$ be an integer such that $F(x)$ is $(x-s)$-regular with respect to $3$ and $s\equiv -1\md{3}$, then $F(x)=220s^9(x-s)^3+66s^{10}(x-s)^2+(a+12s^{11})(x-s)+b+as+s^{12}$. Since $\nu_3(66s^{10})=1$, then $3$ divides $i(K)$ if and only if  $\nu_3(b+as+s^{12})>2\nu_3(a+12s^{11})-1$ or $\nu_3(b+as+s^{12})<2\nu_3(a+12s^{11})-1$, $\nu_3(b+as+s^{12})$ is odd, and $(b+as+s^{12})_3\equiv -1\md{3}$. In this case, $3\mathbb{Z}_{K}=\mathfrak{p}_1^3\mathfrak{p}_{211}\mathfrak{p}_{212}\mathfrak{p}_{22}\mathfrak{p}_3^3$ with $f_1=f_{211}=f_{212}=f_{22}=1$ and $f_3=2$. Hence $3$ divides $i(K)$ and by \cite[Corollary, page 230]{En}, $\nu_3(i(K))=1$. 
					
					\item[(e)] If $(a,b)\in\{(3,29),(3,56),(30,2),(30,56),(57,2),(57,29)\}\md{81}$, then $N_{\phi_2}^+(F)=S_{2}$ has a single side joining $(0,3)$ and $(3,0)$ with $R_{\lambda_2}(F)(y)=-y^3+y^2-y-1$ or $R_{\lambda_2}(F)(y)=-y^3+y^2-y+1=-(y-1)(y^2+1)$. Thus
					$3\mathbb{Z}_{K}=\mathfrak{p}_1^3\mathfrak{p}_2\mathfrak{p}_3^2$ with $f_1=1$, $f_2=3$, and $f_3=2$ or $3\mathbb{Z}_{K}=\mathfrak{p}_1^3\mathfrak{p}_{21}\mathfrak{p}_{22}\mathfrak{p}_3^2$ with $f_1=f_{21}=1$ and $f_{22}=f_3=2$ respectively. Hence $3$ does not divide $i(K)$.
					
					\item[(f)] If $(a,b)\in\{(3,2),(30,29),(57,56)\}\md{81}$, then $N_{\phi_2}^+(F)=S_{21}+S_{22}$ has two sides joining $(0,\mu)$, $(2,2)$, and $(3,0)$ with $\mu\geq4$ and $R_{\lambda_{22}}(F)(y)=-y^2+y-1=-(y+1)^2\in\mathbb{F}_{\phi_2[y]}$. Let $s$ be an integer such that $F(x)$ is $(x-s)$-regular with respect to $3$ and $s\equiv -1\md{3}$, then $F(x)=220s^9(x-s)^3+66s^{10}(x-s)^2+(a+12s^{11})(x-s)+b+as+s^{12}$. Since $\nu_3(66s^{10})=1$, then $3$ divides  $i({K})$ if and only if  $\nu_3(b+as+s^{12})>2\nu_3(a+12s^{11})-1$ or $\nu_3(b+as+s^{12})<2\nu_3(a+12s^{11})-1$, $\nu_3(b+as+s^{12})$ is odd, and $(b+as+s^{12})_3\equiv -1\md{3}$. In this case, $3\mathbb{Z}_{K}=\mathfrak{p}_1^3\mathfrak{p}_{211}\mathfrak{p}_{212}\mathfrak{p}_{22}\mathfrak{p}_3^3$ with $f_1=f_{211}=f_{212}=f_{22}=1$ and $f_3=2$. Hence $3$ divides $i(K)$ and by \cite[Corollary, page 230]{En}, $\nu_3(i(K))=1$. 
					
					\item[(g)] If $(a,b)\in\{(21,47),(48,74),(75,20)\}\md{81}$, then $N_{\phi_2}^+(F)=S_2$ has a single side joining $(0,3)$ and $(3,0)$ with $R_{\lambda_2}(F)(y)=-y^3+y^2+y+1$, which is irreducible over $\mathbb{F}_{\phi_2}$. Thus $3\mathbb{Z}_{K}=\mathfrak{p}_1^3\mathfrak{p}_2\mathfrak{p}_3^3$ with $f_1=1$, $f_2=3$, and $f_3=2$. Hence $3$ does not divide $i(K)$. 
					
					\item [(h)] If $(a,b)\in\{(21,20),(48,47),(75,74)\}\md{81}$, then $N_{\phi_2}^+(F)=S_{21}+S_{22}$ has two sides joining $(0,\mu)$, $(1,2)$, and $(3,0)$ with $\mu\geq4$ and $R_{\lambda_{22}}(F)(y)=-y^2+y+1$, which is irreducible over $\mathbb{F}_{\phi_2}$. Thus $3\mathbb{Z}_{K}=\mathfrak{p}_1^3\mathfrak{p}_{21}\mathfrak{p}_{22}\mathfrak{p}_{3}^3$ with $f_1=f_{21}=1$ and $f_{22}=f_3=2$.
					
					\item[(i)] If $(a,b)\in\{(21,74),(48,20),(75,47)\}\md{81}$, then $N_{\phi_2}^+(F)=S_2$ has a single side joining $(0,3)$ and $(3,0)$ with $R_{\lambda_2}(F)(y)=-y^3+y^2+y-1=-(y+1)(y-1)^2\in\mathbb{F}_{\phi_2}[y]$. Let $s$ be an integer such that $F(x)$ is $(x-s)$-regular with respect to $3$ and $s\equiv -1\md{3}$, then $F(x)=220s^9(x-s)^3+66s^{10}(x-s)^2+(a+12s^{11})(x-s)+b+as+s^{12}$. Since $\nu_3(66s^{10})=1$, then $3$ divide $i(K)$ if and only if  $\nu_3(b+as+s^{12})>2\nu_3(a+12s^{11})-1$ or $\nu_3(b+as+s^{12})<2\nu_3(a+12s^{11})-1$, $\nu_3(b+as+s^{12})$ is odd, and $(b+as+s^{12})_3\equiv -1\md{3}$. In this case, $3\mathbb{Z}_{K}=\mathfrak{p}_1^3\mathfrak{p}_{211}\mathfrak{p}_{212}\mathfrak{p}_{22}\mathfrak{p}_3^3$ with $f_1=f_{211}=f_{212}=f_{22}=1$ and $f_3=2$. Hence $3$ divides $i(K)$ and by \cite[Corollary, page 230]{En}, $\nu_3(i(K))=1$. 
					
				\end{enumerate}				
			\end{enumerate}
			\item If $\nu_3(b)=0$ and $b\equiv 1\md{3}$, then $F(x)\equiv (x^2+x-1)^3(x^2-x-1)^3\md{3}$. Let $\phi_1=x^2+x-1$ and $\phi_2=x^2-x-1$, then 
			$$
			\begin{array}{lll}
				F(x)&=&\cdots+(338-256x)\phi_1^3+(468-474x)\phi_1^2+(324-420x)\phi_1+(a-144)x+b+89,\\
				&=&\cdots+(338+256x)\phi_2^3+(468+474x)\phi_2^2+(324+420x)\phi_2+(a+144)x+b+89.
			\end{array}
			$$
			\begin{enumerate}
				\item[(i)] If $(a,b)\equiv (0,1)\md{9}$, then by Theorem $\ref{thm1}$, $3$ divides $i(K)$.
				
				\item[(ii)] If $(a,b)\not\equiv (0,1)\md{9}$, then by Theorem $\ref{mong}$, $3$ does not divide $i(K)$.
			\end{enumerate}
			
			\item If $\nu_3(b)\geq 1$, then $F(x)\equiv x^{12}\md{3}$. Let $\phi=x$, then $F(x)=\phi^{12}+a\phi+b$.\\
			\textbf{If $\nu_3(b)>\cfrac{12}{11}\nu_3(a)$}, then $N_{\phi}(F)$ has two sides of degree $1$ each, because $\nu_3(a)\leq 10$ by assumption. Thus $3\mathbb{Z}_{K}=\mathfrak{p}_1\mathfrak{p}_{22}^11$ with residue degree $1$ each. Hence $3$ does not divide $i(K)$.\\
			\textbf{If $\nu_3(b)\geq\cfrac{12}{11}\nu_3(a)$}, then we have the following cases:
			\begin{enumerate}
				\item[(i)] If gcd$(\nu_3(b),12)=1$, then $N_{\phi}(F)$ has a single side of degree $1$. Thus $3\mathbb{Z}_{K}=\mathfrak{p}^{12}$ with residue degree $1$. Hence $3$ does not divide $i(K)$.
				
				\item[(ii)] If $\nu_3(b)\in\{2,10\}$, then $N_{\phi}(F)=S$ has a single side  joining $(0,\nu_3(b))$ and $(12,0)$ with $R_{\lambda}(F)(y)=y^2\pm1\in\mathbb{F}_{\phi}[y]$. Thus the possible cases are: $3\mathbb{Z}_{K}=\mathfrak{p}^6$ with $f=2$ and $3\mathbb{Z}_{K}=\mathfrak{p}_1^6\mathfrak{p}_2^6$ with residue degree $1$ each. Hence $3$ does not divide $i(K)$.
				
				\item[(iii)] If $\nu_3(b)\in\{3,9\}$, then $N_{\phi}(F)=S$ has a single side  joining $(0,\nu_3(b))$ and $(12,0)$ with $R_{\lambda}(F)(y)=(y+b_3)^3\in\mathbb{F}_{\phi}[y]$. Since $-\lambda=1/4$ is the slope of $S$, then $4$ divides the ramification index $e$ of any prime ideal of $\mathbb{Z}_{K}$ lying above $3$. Thus the possible cases are: $3\mathbb{Z}_{K}=\mathfrak{p}^4$ with $f=3$, $3\mathbb{Z}_{K}=\mathfrak{p}_1^4\mathfrak{p}_2^4$ with $f_1=2$ and $f_2=1$, $3\mathbb{Z}_{K}=\mathfrak{p}^{12}$ with $f=1$, $3\mathbb{Z}_{K}=\mathfrak{p}_1^8\mathfrak{p}_2^4$ with $f_1=f_2=1$, and $3\mathbb{Z}_{K}=\mathfrak{p}_1^4\mathfrak{p}_2^4\mathfrak{p}_3^4$ with $f_1=f_2=f_3=1$. Hence, $3$ does not divide $i(K)$.
				
				\item[(iv)] If $\nu_3(b)=4$, then $N_{\phi}(F)=S$ has a single side joining $(0,4)$ and $(12,0)$ with $R_{\lambda}(F)(y)=y^4+b_3\in\mathbb{F}_{\phi}[y]$. Thus $R_{\lambda}(F)(y)=(y+1)(y-1)(y^2+1)$ or $R_{\lambda}(F)(y)=(y^2+y-1)(y^2-y-1)$, and so $3\mathbb{Z}_{K}=\mathfrak{p}_1^3\mathfrak{p}_2^3\mathfrak{p}_3^3$ with $f_1=3$ and $f_2=f_3=1$ or $3\mathbb{Z}_{K}=\mathfrak{p}_1^3\mathfrak{p}_2^3$ with $f_1=f_2=3$ respectively. Hence $3$ does not divide $i(K)$.
				
				\item[(v)] If $\nu_3(b)=6$; $b\equiv \pm3^6\md{3^7}$ and $a\equiv 0\md{3^6}$, then $N_{\phi}(F)=S$ has a single side joining $(0,6)$ and $(12,0)$ with $R_{\lambda}(F)(y)=y^6+b_3\in\mathbb{F}_{\phi}[y]$.\\
				For $b_3\equiv1\md{3}$; $b\equiv 3^6\md{3^7}$, we have $R_{\lambda}(F)(y)=(y^2+1)^3$. Since deg$(R_{\lambda}(F))=2$ and $-\lambda=-1/2$ is the slope of $S$, then $2$ divides the ramification index $e$ and the residue degree $f$ of any prime ideal of $\mathbb{Z}_{K}$ lying above $3$. Thus the possible cases are: $3\mathbb{Z}_{K}=\mathfrak{p}^6$ with $f=2$,  $3\mathbb{Z}_{K}=\mathfrak{p}_1^4\mathfrak{p}_2^2$ with $f_1=f_2=2$, $3\mathbb{Z}_{K}=\mathfrak{p}_1^2\mathfrak{p}_2^2\mathfrak{p}_3^2$ with $f_1=f_2=f_3=2$, $3\mathbb{Z}_{K}=\mathfrak{p}_1^2\mathfrak{p}_2^2$ with $f_1=4$ and $f_2=2$, and $3\mathbb{Z}_{K}=\mathfrak{p}^2$ with $f=6$. Hence $3$ does not divide $i(K)$.\\
				For $b_3\equiv -1\md{3}$; $b\equiv -3^6\md{3^7}$, we have
				$R_{\lambda}(F)(y)=(y-1)^3(y+1)^3$. In this case, we have to use Newton polygon of second order. Let $\omega_2=2[\nu_3,\phi,1/2]$ be the valuation of second order Newton polygon, $g_{2,1}=x^2-3$, and $g_{2,2}=x^2+3$ be the key polynomials of $(\omega_2,y-1)$ and $(\omega_2,y+1)$ respectively, then
				$$
				\begin{array}{lll}
					F(x)&=&\cdots+540g_{2,1}^3+1215g_{2,1}^2+1458g_{2,1}+ax+b+729,\\
					&=&\cdots-540g_{2,2}^3+1215g_{2,2}^2-1458g_{2,2}+ax+b+729.
				\end{array}
				$$
				We have $\omega_2(x)=1$, $\omega_2(g_{2,1})=\omega_2(g_{2,2})=2$, and $\omega_2(m)=2\times \nu_3(m)$ for every $m\in\mathbb{Q}_3$.
				\begin{enumerate}
					\item[(a)] If $a\equiv \pm3^6\md{3^7}$, then $N_{2,i}^+(F)=T_{i}$ has a single side  joining $(0,13)$ and $(3,12)$ for every $i=1,2$. Thus $3\mathbb{Z}_{K}=\mathfrak{p}_1^6\mathfrak{p}_2^6$ with residue degree $1$ each. Hence $3$ does not divide $i(K)$.
					
					\item[(b)] If $b\equiv -3^6\pm3^7\md{3^8}$ and $a\equiv 0\md{3^7}$, then $N_{2,i}^+(F)=T_{i}$ has a single side  joining $(0,14)$ and $(3,12)$ for every $i=1,2$. Thus $3\mathbb{Z}_{K}=\mathfrak{p}_1^6\mathfrak{p}_2^6$ with residue degree $1$ each. Hence $3$ does not divide $i(K)$.
					
					\item[(c)] If $b\equiv -3^6\md{3^8}$ and $a\equiv \pm3^7\md{3^8}$, then $N_{2,i}^+(F)=T_{i}$ has a single side  joining $(0,15)$ and $(3,12)$ with $R_{2i}(y)=-y^3-y^2\pm1=-(y\pm1)(y^2\mp y-1)\in\mathbb{F}_{2i}[y]$. Thus $3\mathbb{Z}_{K}=\mathfrak{p}_{11}^2\mathfrak{p}_{12}^2\mathfrak{p}_{21}^2\mathfrak{p}_{22}^2$ with $f_{11}=f_{21}=1$ and $f_{12}=f_{22}=2$. Hence $3$ does not divide $i(K)$.
					
					\item[(d)] If $b\equiv -3^6 \md{3^8}$ and $a\equiv 0\md{3^8}$, then $N_{2,i}^+(F)=T_{i1}+T_{i2}$ has two sides joining $(0,u)$, $(1,14)$ with $u\geq 16$, and $(3,12)$ with $d(T_{1i})=1$, $R_{212}(y)=-y^2-1$, and $R_{222}(y)=y^2+1$, which are irreducible over $\mathbb{F}_{21}$ and $\mathbb{F}_{22}$ respectively. Thus $3\mathbb{Z}_{K}=\mathfrak{p}_{11}^2\mathfrak{p}_{12}^2\mathfrak{p}_{21}^2\mathfrak{p}_{22}^2$ with $f_{11}=f_{21}=1$ and $f_{12}=f_{22}=2$. Hence $3$ does not divide $i(K)$.
				\end{enumerate}
			\end{enumerate}
		\end{enumerate}
		\begin{flushright}
			$\square$
		\end{flushright}
		
		\textit{\textbf{Proof of Theorem \ref{pro1}}}.\\
		For $p=5$. Since $\Delta(F)=2^{24}\times 3^{12}b^{11}-11^{11}a^{12}$ and thanks to the index formula $(1.1)$, if $(a,b)\notin\{(0,0),(1,1),(2,1),(3,1),(4,1)\}\md{5}$, then $5$ does not divide $(\mathbb{Z}_{K}:\mathbb{Z}[\alpha])$, and so $5$ does not divide $i(K)$. Assume that $(a,b)\in\{(0,0),(1,1),(2,1),(3,1),(4,1)\}\md{5}$, then we have the following cases:
		\begin{enumerate}
			\item If $(a,b)\in\{(1,1),(2,1),(3,1),(4,1)\}\md{5}$, then $F(x)\equiv \phi_{i1}\cdot\phi_{i2}\cdot\phi_{i3}^2\md{5}$ with deg$(\phi_{i1})=6$, deg$(\phi_{i2})=4$, deg$(\phi_{i3})=1$, and $\phi_{ij}$ is irreducible over $\mathbb{F}_5$ for $i=1,\dots,4$ and $j=1,2,3$. Thus the possible cases are: $5\mathbb{Z}_{K}=\mathfrak{p}_1\mathfrak{p}_2\mathfrak{p}_3$ with $f_1=6$, $f_2=4$, and $f_3=2$, $5\mathbb{Z}_{K}=\mathfrak{p}_1\mathfrak{p}_2\mathfrak{p}_3^2$ with $f_1=6$, $f_2=4$, and $f_3=1$, and $5\mathbb{Z}_{K}=\mathfrak{p}_1\mathfrak{p}_2\mathfrak{p}_{31}\mathfrak{p}_{32}$ with $f_1=6$, $f_2=4$, and $f_{31}=f_{32}=1$. Hence $3$ does not divide $i(K)$.
			\item If $(a,b)\equiv (0,0)\md{5}$, then $F(x)\equiv x^{12}\md{5}$. Let $\phi=x$, then $F(x)=\phi^{12}+a\phi+b$.\\
			\textbf{If $\nu_5(b)>\cfrac{12}{11}\nu_5(a)$}, then $N_{\phi}(F)$ has two sides of degree $1$ each, because $\nu_5(a)\leq10$ by assumption. Thus $5\mathbb{Z}_{K}=\mathfrak{p}_1\mathfrak{p}_2^{11}$ with residue degree $1$ each. Hence $5$ does not divide $i(K)$.\\
			\textbf{If $\nu_5(b)\geq\cfrac{12}{11}\nu_5(a)$}, then we have the following cases:
			\begin{enumerate}
				\item[(a)] If $\gcd(\nu_5(b),12)=1$, then $N_{\phi}(F)$ has a single side of degree $1$. Thus $5\mathbb{Z}_{K}=\mathfrak{p}^{12}$ with residue degree $1$. Hence $5$ does not divide $i(K)$.
				
				\item[(b)] If $\nu_5(b)\in\{2,3,4,8,9,10\}$, then $N_{\phi}(F)=S$ has a single side joining $(0,\nu_3(b))$ and $(12,0)$ with slopes $-\lambda=-1/6,-1/4,-1/3,-4/3,-3/4,-5/6$ respectively. Since $e_i\in\{3,4,6\}$ divides the ramification index $e$ of any prime ideal of $\mathbb{Z}_{K}$ lying above $5$, then $\phi$ can provide at most four prime ideals of $\mathbb{Z}_{K}$ lying above $5$ with residue degree $1$. Since $\mathcal{N}_1(5)<\mathcal{N}_f(5)$ for every integer $f\geq2$, then $5$ does not divide $i(K)$.
				
				\item[(c)] If $\nu_5(b)=6$, then $N_{\phi}(F)=S$ has a single side joining $(0,6)$ and $(12,0)$ with $R_{\lambda}(F)(y)=y^6+b_5$, which is separable over but does not split completely over $\mathbb{F}_{\phi}$. In particular, $\phi$ can provide at most two prime ideal of $\mathbb{Z}_{K}$ lying above $5$ with residue degree $1$. Hence, $5$ does not divide $i(K)$.
			\end{enumerate}
		\end{enumerate}
		
		For $p=7$. Since $\Delta(F)=2^{24}\times 3^{12}b^{11}-11^{11}a^{12}$ and thanks to the index formula $(1.1)$, if $(a,b)\notin\{(0,0),(1,4),(2,4),(3,4),(4,4),(5,4),(6,4)\}\md{7}$, then $7$ does not divide $(\mathbb{Z}_{K}:\mathbb{Z}[\alpha])$, and so $7$ does not divide $i(K)$. Assume that $(a,b)\in\{(0,0),(1,4),(2,4),(3,4),$\\
		$(4,4),(5,4),(6,4)\}\md{7}$, then we have the following cases:
		
		\begin{enumerate}
			\item If $(a,b)\in\{(1,4),(2,4),(3,4),(4,4),(5,4),(6,4)\}\md{7}$  then $F(x)\equiv \phi_{i1}\cdot\phi_{i2}\cdot\phi_{i3}^2\md{7}$ with deg$(\phi_{i1})=7$, deg$(\phi_{i2})=3$, deg$(\phi_{i3})=1$, and $\phi_{ij}$ is irreducible over $\mathbb{F}_{7}$ for $i=1,\dots,6$ and $j=1,2,3$. Thus  the possible cases are: $7\mathbb{Z}_{K}=\mathfrak{p}_1\mathfrak{p}_2\mathfrak{p}_3$ with $f_1=7$, $f_2=3$, and $f_3=2$, $7\mathbb{Z}_{K}=\mathfrak{p}_1\mathfrak{p}_2\mathfrak{p}_3^2$ with $f_1=7$, $f_2=3$, and $f_3=1$, and $7\mathbb{Z}_{K}=\mathfrak{p}_1\mathfrak{p}_2\mathfrak{p}_{31}\mathfrak{p}_{32}$ with $f_1=7$, $f_2=3$, and $f_{31}=f_{32}=1$. Hence $7$ does not divide $i(K)$.
			
			\item If $(a,b)\equiv (0,0)\md{7}$, then $F(x)\equiv x^{12}\md{7}$. Let $\phi=x$, then $F(x)=\phi^{12}+a\phi+b$.
			\begin{enumerate}
				\item[(a)] If $\nu_7(b)>\cfrac{12}{11}\nu_7(a)$, then $N_{\phi}(F)$ has two sides of degree $1$ each, because $\nu_7(a)\leq 10$ by assumption. Thus $7\mathbb{Z}_{K}=\mathfrak{p}_1\mathfrak{p}_{2}^{11}$. Hence $7$ does not divide $i(K)$.
				
				\item[(b)] If $\nu_7(b)>\cfrac{12}{11}\nu_7(a)$, then $N_{\phi}(F)$ has a single side with slope $-\lambda_i=-h_i/e_i$ such that $e_i\in\{2,3,4,6,12\}$ and gcd$(h_i,e_i)=1$. Since $e_i$ divides the ramification index $e$ of any prime ideal of $\mathbb{Z}_{K}$ lying above $7$, then $\phi$ can provide at most six prime ideals of $\mathbb{Z}_{K}$ lying above $7$ with residue degree $1$. Also, we have $\mathcal{N}_1(7)<\mathcal{N}_f(7)$ for any integer $f\geq2$. We conclude that $7$ does not divide $i(K)$. 
			\end{enumerate}
		\end{enumerate}
		
		For $p=11$. Since $\Delta(F)=2^{24}\times 3^{12}b^{11}-11^{11}a^{12}$ and thanks to the index formula $(1.1)$, if $11$ does not divide $b$, then $11$ does not divide $(\mathbb{Z}_{K}:\mathbb{Z}[\alpha])$, and so $7$ does not divide $i(K)$. Assume that $11$ divides $b$, then we have the following cases:
		\begin{enumerate}
			\item If $\nu_{11}(a)=0$, then $F(x)\equiv x(x+a)^{11}\md{11}$. Let $\phi_1=x$ and $\phi_2=x+a$, then $\phi_1$ provides a unique prime ideal $\mathfrak{p}_1$ of $\mathbb{Z}_{K}$ lying above $11$ with residue degree $1$. Let $F(x)=\phi_2^{12}-12a\phi_2^{11}+66a^2\phi_2^{10}-220a^3\phi_2^9+495a^4\phi_2^8-792a^5\phi_2^7+924a^6\phi_2^6-792a^7\phi_2^5+495a^8\phi_2^4-220a^9\phi_2^3+66a^{10}\phi_2^2+(a-12a^{11})\phi_2+b-a^2+a^{12}$.
			\begin{enumerate}
				\item[(a)] If $\nu_{11}(b-a^2+a^{12})=1$, then $N_{\phi_2}^+(F)$ has a single side of height $1$. Thus $11\mathbb{Z}_{K}=\mathfrak{p}_1\mathfrak{p}_2^{11}$ with residue degree $1$ each. Hence $11$ does not divide $i(K)$.
				
				\item[(b)] If $\nu_{11}(b-a^2+a^{12})\geq2$ and $\nu_{11}(a-12a^{11})=1$, then $N_{\phi_2}^+(F)$ has two sides of degree $1$ each. Thus $11\mathbb{Z}_{K}=\mathfrak{p}_1\mathfrak{p}_{21}\mathfrak{p}_{22}^{10}$ with residue degree $1$ each. Hence $11$ does not divide $i(K)$.
				
				\item[(c)] For $\nu_{11}(b-a^2+a^{12})\geq2$ and $\nu_{11}(a-12a^{11})\geq2$. Since $\nu_{11}(66a^{10})=1$, then $\phi_2$ provides a prime ideal of $\mathbb{Z}_{K}$ lying above $11$ with residue degree $1$ and ramification index $9$. Thus $\phi_2$ can provide at most three prime ideals of $\mathbb{Z}_K$ lying above $11$ with residue degree $1$ each. Since $\mathcal{N}_{1}(11)<\mathcal{N}_f(11)$ for every integer $f\geq2$, we conclude that $11$ does not divide $i(K)$.
			\end{enumerate}
			
			\item If $\nu_{11}(a)\geq 1$, then $F(x)\equiv x^{12}\md{11}$. Let $\phi=x$, then $F(x)=\phi^{12}+a\phi+b$.
			\begin{enumerate}
				\item[(a)] If $\nu_{11}(b)>\cfrac{12}{11}\nu_{11}(a)$, then $N_{\phi}(F)$ has two sides of degree $1$ each, because $\nu_{11}(a)\leq 10$ by assumption. Thus $11\mathbb{Z}_{K}=\mathfrak{p}_1\mathfrak{p}_2^{11}$ with residue degree $1$ each. Hence $11$ does not divide $i(K)$.
				
				\item[(b)] If $\nu_{11}(b)<\cfrac{12}{11}\nu_{11}(a)$, then $N_{\phi}(F)$ has a single side with slope $-\lambda_i=-h_i/e_i$ such that $e_i\in\{2,3,4,6,12\}$ and gcd$(h_i,e_i)=1$. Since $e_i$ divides the ramification index $e$ of any prime ideal of $\mathbb{Z}_{K}$ lying above $11$, then $\phi$ can provide at most six prime ideals of $\mathbb{Z}_{K}$ lying above $11$ with residue degree $1$. Also, we have $\mathcal{N}_1(11)<\mathcal{N}_f(11)$ for any integer $f\geq2$. In this case, $11$ does not divide $i(K)$. 
			\end{enumerate}
			
		\end{enumerate}
		For $p\geq13$, since there are at most $12$ prime ideals of $\mathbb{Z}_{K}$ lying above $p$ with residue degree $1$ each, and there are at least $p\geq 13$ monic irreducible polynomial of degree $f$ in $\mathbb{F}_p[x]$ for every positive integer $f$; we conclude that $p$ does not divide $i(K)$. 
		\begin{flushright}
			$\square$
		\end{flushright}
		\section{Examples}
		Let $F(x)=x^{12}+ax^m+b\in \mathbb{Z}[x]$ be a monic irreducible polynomial and $K$ the number field generated by a  root of $F(x)$.
		\begin{enumerate}
			\item For $a=72$ and $b=51$, we have $(a,b)\equiv(0,3)\md{8}$, then by Theorem $\ref{thm1}$, $K$ is not monogenic for every $m=1,\dots,11$.
			
			\item For $a=108$ and $b=26$, we have $(a,b)\equiv(0,-1)\md{9}$, then by Theorem $\ref{thm1}$, $K$ is not monogenic for every $m=1,\dots,11$.
			
			\item For $a=60$, $b=15$, and $m=6$, we have $(a,b)\equiv(4,7)\md{8}$, then by Theorem $\ref{thm2}$, $K$ is not monogenic.
			
			\item For $a=14$, $b=90$, and $m=10$, we have $a\equiv -1\md{3}$, $\nu_3(b)=2$, and $b_3\equiv 1\md{3}$, then by Corollary $\ref{cor}$, $K$ is not monogenic.

			\item For $a=576$, $b=386$, and $m=1$, we have  $(a,b)\equiv(64,112)\md{256}$, then by Theorem $\ref{thmp2}$, $2$ divides $i(K)$, and so $K$ is not monogenic.
			
			\item For $a=84$, $b=147$, and $m=1$, we have $F(x)\equiv (x-1)^4(x^2+x+1)^2\md{2}$. Let $F(x)=\cdots+(x-1)^4+(x-1)^3+(x-1)^2+(x-1)+232$, then $N_{x-1}^+(F)=S_1+S_2$ has two sides joining $(0,3)$, $(2,1)$, and $(4,0)$ with $R_{\lambda_1}(F)(y)=(y+1)^2\in\mathbb{F}_{(x-1)}[y]$, and so $F(x)$ is not $2$-regular. The regular integer in this case is $s=3$. Let $F(x)=\cdots+3247695(x-3)^4+4330260(x-3)^3+3897234(x-3)^2+2125848(x-3)+531840$, then $N_{x-3}^+(F)=S_1+S_2+S_3$ has three sides of degree $1$ each. Thus $(x-3)$ provides three prime ideals of $\mathbb{Z}_{K}$ lying above $2$ with residue degree $1$ each, and so $2$ divides $i(K)$. By Theorem \ref{thmp2} (2), $\nu_2(i(K))=1$. On the other hand, $\nu_3(a)\geq1$ and $\nu_3(b)=1$, then $N_x(F)$ has a single side of height $1$. Thus $3$ does not divide $(\mathbb{Z}_K:\mathbb{Z}[\alpha])$. We conclude that $i(K)=2$ and so $K$ is not monogenic.
			
			\item For $a=120$, $b=-610$, and $m=1$, we have  $(a,b)\equiv(39,38)\md{81}$. Since $u=\nu_3(b-a+1)=\nu_3(-729)=6$ and $v=\nu_3(a-12)=\nu_3(108)=3$, by Theorem $\ref{thmp3}$, $3$ is a divides $i(K)$ with $\nu_3(i(K))=1$.  On the other hand, $\nu_2(a)\geq1$ and $\nu_2(b)=1$, then $N_x(F)$ has a single side of height $1$. Thus $2$ does not divide $(\mathbb{Z}_K:\mathbb{Z}[\alpha])$. We conclude that $i(K)=3$ and so $K$ is not monogenic.
			
			\item For $a=45996$, $b=373907$, and $m=1$, we have  $(a,b)\equiv (4,3)\md{8}$. Since $u=\nu_2(a+b+1)=\nu_2(419904)=6$ and $v=\nu_2(a+12)=\nu_2(46008)=3$, by Theorem \ref{thmp2} (2), $2$ divides $i(K)$ with $\nu_2(i(K))=1$. 0n the other hand, we have $(a,b)\equiv(69,11)\md{81}$, $u=\nu_3(a+b+1)=\nu_3(419904)=8$, and $v=\nu_3(a+12)=\nu_3(46008)=4$, then by Theorem $\ref{thmp3}$ (2), $3$ divides $i(K)$ with $\nu_3(i(K))=1$. We conclude that $i(K)=6$ and so $K$ is not monogenic.
			
			\item For $a=336$, $b=5195$, and $m=1$, we have $(a,b)\equiv (0,3)\md{8}$, then by Theorem $\ref{thmp2}$ (1), $2$ divides $i(K)$ with $\nu_2(i(K))=2$. On the other hand, $(a,b)\equiv(12,11)\md{81}$, $\mu=\nu_3(a-b+1)=\nu_3(4860)=5$, $\tau=\nu_3(a-12)=\nu_3(324)=4$, and $(a-b+1)_3=(4860)_3\equiv -1\md{3}$, then by Theorem $\ref{thmp3}$ (4), $3$ divides $i(K)$ with $\nu_3(i(K))=1$. We conclude that $i(K)=12$ and so $K$ is not monogenic.
			
			\item For $a=24576$, $b=216128$, and $m=1$, we have $(a,b)\equiv (0,64)\md{1024}$, then by Theorem $\ref{thmp2}$ (4), $2$ divides $i(K)$ with $\nu_2(i(K))=3$. On the other hand, $(a,b)\equiv (42,38)\md{81}$, then  by Theorem $\ref{thmp3}$ (3), $F(x)$ is not $3$-regular. The regular integer in this case is $s=-2$. Let $F(x)=\cdots-112640(x-s)^3+67584(x-s)^2+171072$. Since $\nu_3(171072)=5$, $\nu_3(67584)=1$, and $(171072)_3\equiv -1\md{3}$, then by Theorem $\ref{thmp3}$ (3), $3$ divides $i(K)$ with $\nu_3(i(K))=1$. We conclude that $i(K)=24$ and so $K$ is not monogenic.

			
			
		\end{enumerate}

	\end{document}